\documentclass{IEEEtran}
\pdfoutput=1

\usepackage{amsthm,amsmath,amssymb,amsfonts}
\usepackage{array}
\let\chapter\section
\usepackage[ruled]{algorithm2e}
\usepackage{enumitem}
\usepackage{graphics,graphicx}
\usepackage{multirow}
\usepackage[square, numbers]{natbib}
\usepackage{subfig}

\usepackage{changepage}

\usepackage{etex}
\reserveinserts{40}

\def\RR{\mathbb{R}}

\def\lone{\ell_1}
\def\lzero{\ell_0}

\def\rfl1{\rho_S^{\lone}(\delta)}

\def\rflqc1{\rho_S^{\lone}(\delta;q,C_1(\delta,\rho)\le\Upsilon)}

\theoremstyle{definition}
\newtheorem{theorem}{Theorem}[section]
\newtheorem{lemma}[theorem]{Lemma}

\newtheorem{definition}[theorem]{Definition}

\def\expanderplain{\mathbb{E}_{k, \varepsilon, d}}
\def\expander{\mathbb{E}_{k, \varepsilon, d}^{m\times n}}
\def\Sparse{\chi_k^n}

\DeclareMathOperator*{\argmax}{arg\,max}
\DeclareMathOperator*{\argmin}{arg\,min}
\DeclareMathOperator{\argsupp}{argsupp}
\DeclareMathOperator{\esp}{\;\;\;\;\;}

\DeclareMathOperator{\median}{median}
\DeclareMathOperator{\mode}{mode}

\DeclareMathOperator{\st}{\mbox{s.t.}}
\DeclareMathOperator{\supp}{supp}

\def\Prob{\Pr}

\def\RR{{\mathbb R}}

\def\supp{\text{supp}}

\def\R{\mathbb{R}}

\def\N{\mathbb{N}}

\newcommand{\nc}{\newcommand}
\nc{\binomial}{\mbox{Bin}}
\nc{\Mod}[1]{\ (\text{mod}\ #1)}

\graphicspath{{figures/}}
\begin{document}
%
\title{Expander $\lzero$-decoding}

\author{Rodrigo Mendoza-Smith
        and~Jared~Tanner
\thanks{Copyright (c) 2015 The Authors. Personal use of this material is permitted. However, permission to use this material for any other purposes must be obtained from the authors.}
\thanks{RMS and JT are with the Mathematical Institute, University of Oxford, Oxford, UK.   RMS is supported by CONACyT. (email: \{mendozasmith,tanner\}@maths.ox.ac.uk)}
\thanks{Manuscript submitted August 2015.}
}

\markboth{Mendoza-Smith and Tanner: \emph{Expander $\lzero$-Decoding}}%
{Mendoza-Smith and Tanner}

\maketitle

\begin{abstract}
We introduce two new algorithms, Serial-$\ell_0$ and Parallel-$\ell_0$ for solving a large underdetermined linear system of equations $y = Ax \in \R^m$ when it is known that $x \in \R^n$ has at most $k < m$ nonzero entries and that $A$ is the adjacency matrix of an unbalanced left $d$-regular expander graph.
The matrices in this class are sparse and allow a highly efficient implementation.   A number of algorithms have been designed to work exclusively under this setting, composing the branch of {\em combinatorial compressed-sensing} (CCS).
Serial-$\ell_0$ and Parallel-$\ell_0$ iteratively minimise $\|y - A\hat x\|_0$ by successfully combining two desirable features of previous CCS algorithms: the information-preserving strategy of ER \citep{jafarpour2008efficient}, and the parallel updating mechanism of SMP \citep{berinde2008practical}. 
We are able to link these elements and guarantee convergence in $\mathcal{O}(dn \log k)$ operations by assuming that the signal is {\em dissociated}, meaning that all of the $2^k$ subset sums of the support of $x$ are pairwise different.
However, we observe empirically that the signal need not be exactly
dissociated in practice.  Moreover, we observe Serial-$\ell_0$ and
Parallel-$\ell_0$ to be able to solve large scale problems with a
larger fraction of nonzeros than other algorithms when the number of
measurements is substantially less than the signal length; in
particular, they are able to reliably solve for a $k$-sparse vector
$x\in\RR^n$ from $m$ expander measurements with $n/m=10^3$ and
$k/m$ up to four times greater than what is achievable by
$\ell_1$-regularization from dense Gaussian measurements.  
Additionally, due to their low computational complexity,
Serial-$\ell_0$ and Parallel-$\ell_0$ are observed to be able to solve
large problems sizes in substantially less time than other algorithms
for compressed sensing. In particular, Parallel-$\ell_0$ is structured to take advantage of massively parallel architectures.
\end{abstract}

\section{Introduction}

Compressed sensing \citep{candes2005decoding, donoho2006compressed, candes2006quantitative, candes2006robust, candes2006stable, candes2006near} considers the problem of sampling and efficiently reconstructing a compressible finite dimensional signal $x \in \R^n$ from far fewer measurements than what Nyquist and Shannon deemed possible \citep{nyquist1928certain, shannon1949communication}.
In its simplest form compressed sensing states that if $x\in\RR^n$ has at most
$k<n$ nonzero entries, then it can be sampled from $m$ linear
measurements $y=Ax\in\RR^m$ and that $x$ can be recovered from $(y,A)$
with computationally efficient algorithms provided $m<n$ is
sufficiently large, see \cite{foucart2013mathematical}.

The most widely studied sensing matrices $A$ are from the classes of:
a) Gaussian or uniformly drawn projections which are most amenable to
precise analysis due to their spherical symmetry, and b) partial Fourier matrices
which have important applications for tomography and have fast
transforms allowing $A$ and $A^*$ to be applied in ${\cal O}(n\log n)$
operations.  Unfortunately the partial Fourier matrices are not known
to allow the asymptotically optimal order number of measurements of
$m\sim k\sim n$, rather the best analysis ensures recovery for $m\sim
k\log^{5}n$ \cite{foucart2013mathematical}.  Sparse binary matrices with a fixed number of
non-zeros per column offer the possibility of $A$ and $A^*$ being
applied in ${\cal O}(n)$ time and for asymptotically optimal order
number of measurements $m\sim k\sim n$ \cite{bah2013vanishingly,berinde2008combining}.
When restricting to these matrices, compressed sensing is referred to as {\em combinatorial
  compressed sensing}, \citep{berinde2008combining}.

\subsection{Combinatorial compressed sensing}

The problem of sparse recovery with compressed sensing resembles the problem of {\em linear sketching} in theoretical computer science.
This area considers {sketching} high dimensional vectors $x \in \R^n$ using a sparse matrix $A \in \R^{m \times n}$ with the aim that $Ax$ has lower dimensionality than $x$, but still preserves some of its properties with high probability.
In an attempt to reconcile this area with the compressed-sensing paradigm, \citep{berinde2008combining} proposed sensing $x \in \Sparse$ using an expander matrix, $i.e.$ the adjacency matrix an unbalanced bipartite graph with high connectivity properties\footnote{See Section \ref{sec:expanders} for details.}.
We denote the $m\times n$ matrices in this class by $\expander$, but abbreviate to $\expanderplain$ when the size is
understood by its context.  Expander matrices $\expander$
are sparse binary matrices with $d \ll m$
ones per column, but with their nonzeros distributed in such a way
that any submatrix composed of $k$ columns has at least
$(1-\varepsilon)kd$ rows which are nonzero \footnote{ Such expander matrices can be generated by drawing i.i.d. columns with the location of their nonzeros drawn uniformly from the $m \choose d$ support sets of cardinality $d$, \cite{bah2013vanishingly}}.
%
This structure makes them suitable for sparse recovery, and also makes them low complexity in terms of storage, generation, and computation (see Table \ref{table:measurement-operators}).
Additionally, some applications like the single-pixel camera
\citep{baraniuk2008single} consider measurement devices with binary
sensors that inherently correspond to binary and sparse inner
products, and that unfortunately, fall outside the set of matrices for
which the widely used restricted isometry techniques apply.

The authors of \citep{berinde2008combining} showed that, although being sparse, expander matrices can sense elements in $\Sparse$ at the optimal measurement rate $\mathcal{O}(k \log(k/m))$, and that these can be recovered accurately and efficiently via $\ell_1$-regularization.
Following this result, a series of algorithms designed specifically to work with expander matrices was presented in \citep{jafarpour2008efficient, berinde2008practical, berinde2009sequential, xu2007efficient}.
The analysis of these algorithms requires the use of techniques and ideas borrowed from combinatorics, which is why this branch was labeled by \citep{berinde2008combining} as {\em combinatorial compressed sensing} (CCS).
It is in this realm that we make our main contributions.

\begin{table}[!htbp]
 \begin{adjustwidth}{-0.5cm}{}
\begin{tabular}{|c|c|c|c|c|}
\hline
& Storage & Generation & $A^*y$& $m$\\
\hline
&&&&\\
{\scriptsize Gaussian/Bernoulli} & $\mathcal{O}(mn)$ & $\mathcal{O}(mn)$ & $\mathcal{O}(mn)$ & $\mathcal{O}(k \log(n/k))$\\
&&&&\\
Partial Fourier & $\mathcal{O}(m)$ & $\mathcal{O}(m)$ & $\mathcal{O}(n\log n)$ &$\mathcal{O}(k \log^5(n))$\\
&&&&\\
Expander & $\mathcal{O}(dn)$ & $\mathcal{O}(dn)$ & $\mathcal{O}(dn)$ &$\mathcal{O}(k \log(n/k))$\\
&&&&\\
\hline
\end{tabular}
\caption{Complexity of  measurement operators.}
\label{table:measurement-operators}
 \end{adjustwidth}
\end{table}

\subsection{Main contributions}
Our work is in the nexus of a series of papers \citep{jafarpour2008efficient, berinde2008practical, berinde2009sequential, xu2007efficient} proposing iterative greedy algorithms for combinatorial compressed sensing. 
The algorithms put forward in the aforementioned sequence of papers
recover the sparsest 
solution of a large underdetermined linear system of equations $y =
Ax$ by iteratively refining an estimation $\hat x$ using information
about the residual $r = y - A\hat x$. 
%
Though these algorithms have the same high-level
perspective\footnote{See Section \ref{sec:prior-art} and Table
  \ref{table:ccs}}, their particulars are optimised to best tradeoff
speed, robustness, and recovery region; see Table~\ref{table:ccs} for
a summary.
For instance, at each iteration, SMP \citep{berinde2008practical} updates several entries of $\hat x$ in parallel, allowing it to provably recover an arbitrary $x \in \Sparse$ in $\mathcal{O}(\log \|x\|_1)$ iterations of complexity $\mathcal{O}(dn + n \log n)$.
However, SMP is only able to
recover the sparsest solution when the fraction of nonzeros in the
signal is substantially less than other compressed sensing algorithms.
%
On the other hand, at each iteration, LDDSR \citep{xu2007efficient} and ER \citep{jafarpour2008efficient} update a single entry of $\hat x$ in such a way that a contraction of $\|y - A\hat x\|_0$ is guaranteed.
This reduction in the residual's sparsity is achieved by exploiting an important property of expander graphs, which we call the {\em information-preserving} property (see Theorem \ref{th:equivalence-lemma}). 
Essentially, this property guarantees that most of the entries from $x$ will appear
repeatedly as entries in $y=Ax$.
%
In other words, it guarantees that for most $i \in [m]$, we will have $y_i \in \{x_j : j \in \supp(x)\}$.
In \citep{xu2007efficient} and \citep{jafarpour2008efficient}, this property is used to give sufficient conditions for decrease of $\|y - A\hat x\|_0$ under the regime of single updating of $\hat x$.  
However, this regime of single updating in LDDSR and ER typically requires greater computational time than existing compressed-sensing algorithms.
Our main contribution is in the design and analysis of an algorithmic model that successfully combines the information-preserving strategy of LDDSR and ER with the parallel updating scheme of SMP.
This synthesis is made possible by assuming that the signal of interest is {\em dissociated}.
%
\begin{definition}[Dissociated signals]
\label{def:dissociated-signals}
A signal $x \in \R^n$ is dissociated if 
\begin{equation}
\label{eq:dissociated}
\sum_{j \in T_1} x_j \neq \sum_{j \in T_2} x_j \esp\forall\; T_1, T_2 \subset \supp(x)\mbox{ $\st$ }T_1 \neq T_2.
\end{equation}
\end{definition}
The name {\em dissociated} comes from the field of additive combinatorics (See Definition 4.32 in \citep{tao2006additive}), where a set $S$ is called {\em dissociated} if the set of all sums of distinct elements of $S$ has maximal cardinality.
Even though the model (\ref{eq:dissociated}) might seem restrictive, it need not be exactly fulfilled for our algorithm to work.
In fact, it is fulfilled almost surely for isotropic signals, and more generally by any signal whose nonzeros can be modelled as being drawn from a continuous distribution.
%
Moreover, it is discussed in Section \ref{sec:dissociated-A} that non-dissociated signals, such as integer or binary signals, can be recovered if instead the columns of $A$ are scaled by dissociated values, and the nonzeros of $x$ are drawn independently of $A$. 
Also, numerical experiments show
that the algorithm recovery ability decreases gracefully as the
dissociated property is lost by having a fraction of the nonzeros in
$x$ be equal, see Figure~\ref{fig:banded_signals},.

With this assumption, our contributions are a form of {\em model-based compressed sensing} \citep{baraniuk2010model} in which apart from
assuming $x \in \Sparse$, 
one also assumes special dependencies between the values of its
nonzeros with the goal to improve the algorithms speed or recovery ability.
Our contributions are Serial-$\ell_0$ and Parallel-$\ell_0$,
Algorithms \ref{alg:serial-l0} and \ref{alg:parallel-l0} respectively,
and their convergence guarantees summarised in Theorem
\ref{th:convergence-expl0de}. 
%
\begin{small}
\begin{algorithm}
 \KwData{$A \in \expander$; $y \in \R^m$; $\alpha \in (1, d]$}
 \KwResult{ $\hat x \in \R^n$ $\st$ $ y = A \hat x$ }
 $\hat x \leftarrow 0$, $r \leftarrow y$\;
 \While{not converged}{
     \For{$j \in [n]$}{
     $T\in \left\{\omega_j \in \R : \|r\|_0-\|r - \omega_j a_j\|_0 \geq \alpha\right\}$\; 
		\For{$\omega_j \in T$}{
			$\hat x_j \leftarrow \hat x_j + \omega_j$\;
		}
		$r \leftarrow y - A\hat x$\;
	}
}
\caption{Serial-$\ell_0$}
\label{alg:serial-l0}
\end{algorithm}
\end{small}
\begin{small}
\begin{algorithm}
 \KwData{$A \in \expander$; $y \in \R^m$; $\alpha \in (1, d]$}
 \KwResult{ $\hat x \in \R^n$ $\st$ $ y = A \hat x$ }
 $\hat x \leftarrow 0$, $r \leftarrow y$\;
  \While{not converged}{
  	$T\leftarrow \left\{(j, \omega_j) \in [n]\times\R : \|r\|_0-\|r - \omega_j a_j\|_0 \geq \alpha\right\}$\; 
	\For{$(j, \omega_j) \in T$}{
		$\hat x_j \leftarrow \hat x_j + \omega_j$\;
	}
	$r\leftarrow y - A\hat x$\;
}
\caption{Parallel-$\ell_0$}
\label{alg:parallel-l0}
\end{algorithm}
\end{small}
\begin{theorem}[Convergence of Expander $\ell_0$-Decoders]
\label{th:convergence-expl0de}
Let $A \in  \expander$ and $\varepsilon \leq 1/4$. and $x \in \Sparse$ be a dissociated signal. Then, Serial-$\ell_0$ and Parallel-$\ell_0$ with $\alpha = (1 - 2\varepsilon)d$ can recover $x$ from $y = Ax \in \R^m$ in $\mathcal{O}(dn\log k)$ operations.
\end{theorem}
The focus of this paper is on charting the development of
Serial-$\ell_0$ and Parallel-$\ell_0$ and on proving Theorem
\ref{th:convergence-expl0de}. 
In doing so, we contrast Serial-$\ell_0$ and Parallel-$\ell_0$ to the
state-of-the-art algorithms for compressed-sensing and show that when
the signal is dissociated, these are the fastest algorithms available
when implemented, respectively, in a serial or a parallel
architecture. 
We support these claims with a series of numerical experiments that
additionally show that any loss in universality due to our signal
model is traded off by unusually high recovery regions when $\delta:=
m/n$ is small and substantially higher than those of previous CCS
algorithms. 
\subsection{Outline}
Section \ref{sec:background} gives the main background theory in expander graphs necessary for our discussion. 
Then, Section \ref{sec:prior-art} reviews past advances in CCS, putting emphasis on deconstructing these into their essential ideas, and on pointing out common elements between them.
Section \ref{sec:contributions} contains our main contributions: Serial-$\ell_0$ and Parallel-$\ell_0$.
We prove Theorem \ref{th:convergence-expl0de} and point out some technical details regarding the implementation of Serial-$\ell_0$ and Parallel-$\ell_0$.
We also discuss some connections of the dissociated model (\ref{eq:dissociated}) with Information Theory.
Finally, in Section \ref{sec:numerical-experiments} we evaluate the empirical performance of these algorithms with a series of numerical experiments. 

\section{Background}
\label{sec:background}
In this section, we present the basic notions of graph theory that are necessary for understanding our subsequent analyses, as well as the relevant previous work in combinatorial compressed sensing. We start by defining some notation.

\subsection{Notation}

For a subset $S \subset \Omega$, we let $|S|$ be its cardinality, and $\Omega \setminus S$ denote its complement.
We adopt notation from combinatorics and use the shorthand $[n] := \{1, \dots, n\}$ for $n\in \mathbb{N}$.
We also define $[n]^{(k)} = \{ S \subset [n] : |S| = k\}$ and $[n]^{(\leq k)} = \{ S \subset [n] : |S| \leq k\}$.
As mentioned in the previous section, for $x\in \R^n$, we let
$\supp(x) = \{i : x_i \neq 0\}$ be its support, and $\argsupp(x) =
\{x_i : i \in \supp(x)\}$ be the set of nonzero values in $x$.
With this, we define $\|x\|_0 = |\supp(x)|$, and $\Sparse = \{ x \in \R^n : \|x\|_0 \leq k\}$; vectors in $\Sparse$ are said to be $k$-sparse.
We let $H_k : \R^n \rightarrow \Sparse$ be the hard thresholding
operator that sets to zero all but the $k$ largest elements in $x$. 
Throughout this work, we implicitly assume that $x \in \R^n$, $y \in \R^m$, and that $A \in \R^{m \times n}$ is a binary sparse matrix with $d$ ones per column.
It is also implicitly assumed that $m < n$ and that $\|x\|_0 < m$. For a given signal $x$, we will use $k$ to refer to its sparsity, unless we specify otherwise.

\subsection{Expander graphs}
\label{sec:expanders}
\begin{figure}[!htbp]
	\centering
	\includegraphics[scale=0.5]{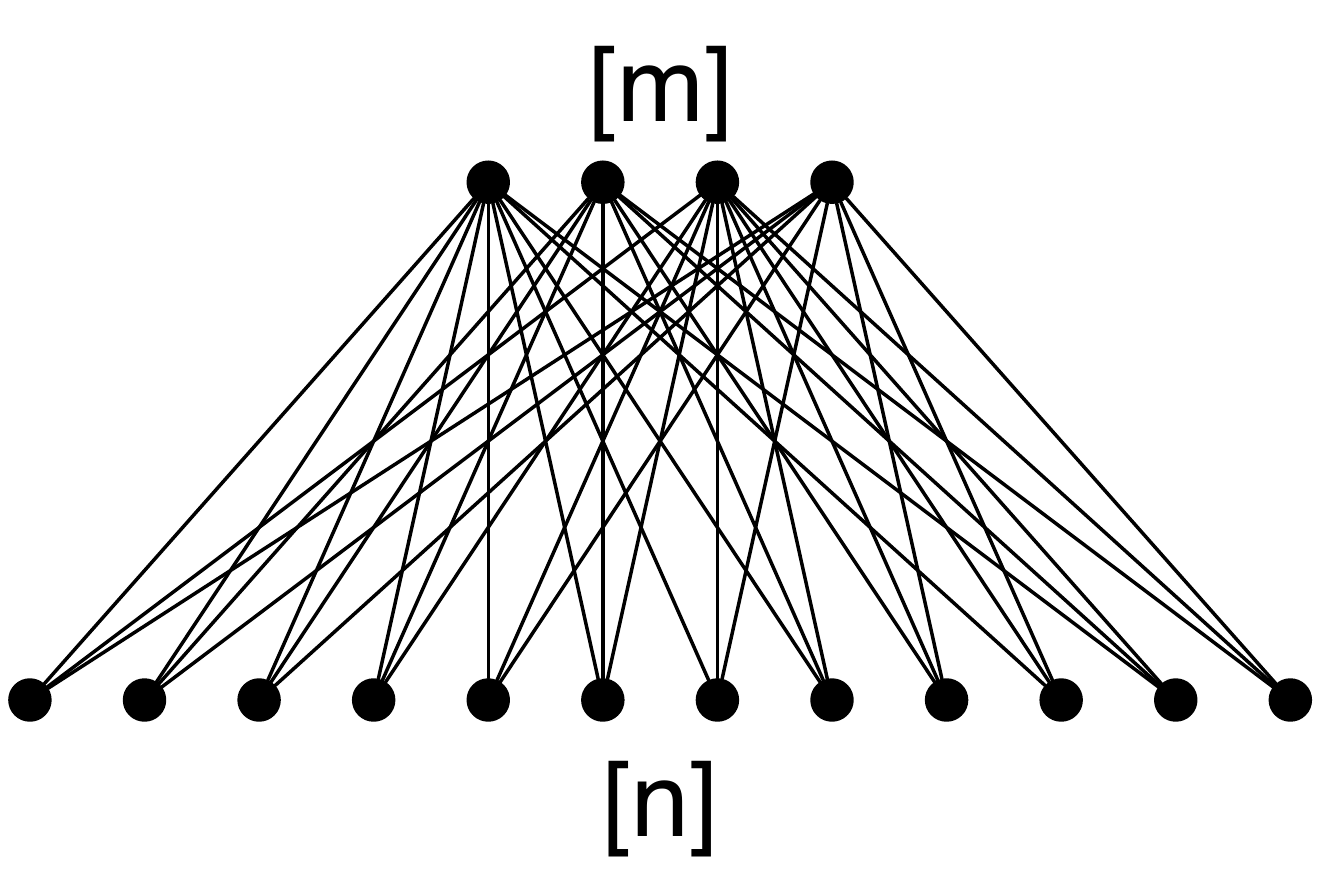}
	\caption[Schematic of expander graph]{{\bf Schematic of an expander graph with $d=3$.} Every left $d$-regular bipartite graph is an expander for some $k$ and $\varepsilon$.}
	\label{fig:schematic_expander}
\end{figure}

A {\em bipartite graph} is a 3-tuple $G = (U,V,E)$ such that $U \cap V = \emptyset$ and $E \subset U \times V$. 
Elements in $U \cup V$ are called {\em nodes}, while tuples in $E$ are called {\em edges}.
Under the assumption that $|U| = n$ and $|V| = m$, we abuse notation and let $U = [n]$ be the set of {\em left-nodes}, and $V = [m]$ be the set of {\em right-nodes}.
A bipartite graph is said to be {\em left $d$-regular} if the number of edges emanating from each left node is identically $d$, and is said to be {\em unbalanced} if $m < n$.
For $S\subset U \cup V$ we define $\mathcal{N}(S) \subset U\cup V$ to be the {\em neighbourhood} of $S$, $i.e.$ the set of nodes in $U\cup V$ that are connected to $S$ through an element of $E$. 
We note that for bipartite graphs, $\mathcal{N}(S) \subset V$ only if $S \subset U$, and $\mathcal{N}(S) \subset U$ only if $S \subset V$.
An {\em expander graph} (Figure \ref{fig:schematic_expander}) is an unbalanced, left $d$-regular, bipartite graph that is {\em well-connected} in the sense of the following definition.
%
\begin{definition}[Expander graph]
\label{def:expander}
An unbalanced, left $d$-regular, bipartite graph $G = ([n],[m],E)$ is a $(k, \varepsilon, d)$-expander if
\begin{equation}
\label{eq:expansion}
|\mathcal{N}(S)| > (1-\varepsilon)d|S| \;\;\forall \;\;S \in [n]^{(\leq k)}.
\end{equation}
We call $\varepsilon \in (0,1)$ the {\em
 expansion parameter} of the graph.
\end{definition}
Hence, the expander graphs that we consider can be thought of as tuples $G = ([n], [m], E)$ such that all subsets $S \in [n]^{(\leq k)}$ have at most $\varepsilon d |S|$ fewer neighbours than the number of edges emanating from $S$. It will be convenient to think of an expander in linear algebra terms, which can be done via its {\em adjacency matrix}.
\begin{definition}[Expander matrix $\expander$]
\label{def:adjacency-expander}
The adjacency matrix of an unbalanced, left $d$-regular, bipartite graph $G = ([n],[m],E)$ is the binary sparse matrix $A \in \R^{m \times n}$ defined as
\begin{equation}
\label{eq:adjacency}
A_{ij} = \left\{ \begin{array}{ll}
1 & \mbox{$i \in \mathcal{N}(j) \subset [m]$}\\
0 & \mbox{ otherwise.}\end{array}\right.
\end{equation}
We let $\expanderplain$ be a the set of adjacency matrices of $(k, \varepsilon, d)$-expander graphs. 
\end{definition}
We note that $A \in \expander$ is a sparse binary matrix with exactly $d$ ones per column, and also that any left $d$-regular bipartite graph will satisfy (\ref{eq:expansion}) for some $k$ and $\varepsilon$.
As mentioned previously, \citep{berinde2008combining} showed that
these matrices possess a bounded restricted isometry constant (RIC) in the
$\ell_1$ norm in the linear growth asymptotic where $k\sim m \sim n\rightarrow \infty$; making these matrices computationally 
highly attractive for compressed sensing.
The existence of expander graphs with optimal measurement rate of $m = \mathcal{O}(k \log(n/k))$, is addressed in the following theorem.
%
\begin{theorem}[Existence of optimal expanders \citep{bassalygo1973complexity, capalbo2002randomness}]
\label{th:expander-existence}
For any $n/2\geq k \geq 1$ and $\varepsilon > 0$, there is a $(k, \varepsilon, d)$-expander with
\begin{equation}
d = \mathcal{O}(\log(n/k)/\varepsilon),\;\;\mbox{ and }\;\; m = \mathcal{O}(k \log(n/k)/\varepsilon^2).
\end{equation}
\end{theorem}
Theorem \ref{th:expander-existence} also implies that in the linear growth
asymptotic of $k\sim m\sim n\rightarrow\infty$ and for a fixed $\varepsilon > 0$, it holds that $d = \mathcal{O}(1)$; that is, the number of nonzeros
per column does not increase with the problem
size.
Apart from this fact, expander matrices are of interest in compressed sensing because they are nearly information preserving, meaning that for $x \in \Sparse$ at least $(1 - 2\varepsilon)kd$ entries of $y = Ax$ equal a nonzero value of $x$. This property is guaranteed by Lemma \ref{th:equivalence-lemma}.
%
%
\begin{lemma}[Information-preserving property]
\label{th:equivalence-lemma}
Let $G = ([n],[m],E)$ be an unbalanced, left $d$-regular, bipartite graph, and $S \in [n]^{(\leq k)}$. Define,
\begin{equation}
\label{eq:N1}
\mathcal{N}_1(S) = \{ i \in \mathcal{N}(S) : |\mathcal{N}(i) \cap S| = 1\},
\end{equation}
and 
\begin{equation}
\label{eq:Ngeq1}
\mathcal{N}_{>1}(S) = \mathcal{N}(S) \setminus \mathcal{N}_1(S).
\end{equation}
Then, $G$ is a $(k, \varepsilon, d)$-expander graph if and only if
\begin{equation}
\label{eq:unique-neighbors}
|\mathcal{N}_1(S)| > (1 - 2\varepsilon)d|S| \;\; \forall \;\; S \in [n]^{(\leq k)}.
\end{equation}
\end{lemma}
\begin{IEEEproof}
See Appendix \ref{app:expander-proofs}.
\end{IEEEproof}

The information-preserving property is widely used in the analysis of
CCS, and is a central piece in the analysis of our algorithms as it
implies the lower $\ell_1$-RIC bound \cite{berinde2008combining}.
Finally, we remark that adjacency matrices of expander graphs are not only useful for compressed-sensing, but also for a number of applications including linear sketching, data-stream computing, graph sketching, combinatorial group testing, network routing, error-correcting codes, fault-tolerance, and distributed storage \citep{berinde2008combining, capalbo2002randomness}.
%

\section{Overview of CCS prior art}
\label{sec:prior-art}
Iterative greedy algorithms for compressed sensing seek the sparsest
solution to a large underdetermined linear system of equations $y =
Ax$ and typically do so by operating on the residual $r = y - A\hat
x$, where $\hat x$ is an estimate of the sparsest solution.
%
Algorithms for combinatorial compressed sensing differ by considering
updating the $j^{th}$ entry of the approximation, $\hat x_j$, based on
a non inner product score $s_j \in \R$ dependent on $r_{\mathcal{N}(j)}$; that is, on the residual restricted to the
support set of the $j^{th}$ column of $A$.
%
In order to standardise the convergence rate guarantees of previous CCS, we
define the notion of an iteration as follows.
%
\begin{definition}[Iteration]
\label{def:iteration}
Let $A \in \R^{m \times n}$, $x \in \R^n$, and $y = Ax$. For an iterative greedy algorithm updating an estimation $\hat x \in \R^n$ of $x$ from a residual $r = y - A\hat x$, an iteration is defined as the sequence of steps performed between two updates of $r$.
\end{definition}
In the remainder of this section we deconstruct past CCS algorithms
into their essential components so as to give a high-level overview of
their shared characteristics. 
\subsection{Sparse Matching Pursuit (SMP)}
SMP was proposed in \citep{berinde2008practical} to decode $\hat x$ from $y=Ax$ with a voting-like mechanism in the spirit of the {\em count-median} algorithm from data-stream computing (see \citep{cormode2005improved} for details).
SMP can also be viewed as an {\em expander} adaptation of the
Iterative Hard Thresholding algorithm (IHT)
\citep{blumensath2009iterative}, which uses the line-search $\hat x
\leftarrow H_k[\hat x + p]$ to minimise $\|y - A\hat x\|_2^2$ over
$\Sparse$, indeed it was rediscovered from this perspective in 
\citep{foucart2013mathematical}[pp. 452] where it is referred to as EIHT.
Due to the structure of expander matrices, SMP chooses the direction $p = \mathcal{M}(y - A\hat x)$ with $\mathcal{M}: \R^m \rightarrow \R^n$ defined as
\begin{equation}
[\mathcal{M}(r)]_j = \median (r_{\mathcal{N}(j)}).
\end{equation}
After thresholding, this choice yields the iteration,
\begin{equation}
\label{eq:smp}
\hat x \leftarrow H_k\left[ \hat x + H_{2k}[\mathcal{M}(y - A\hat x)]\right].
\end{equation}
SMP and its theoretical guarantees are stated in
Algorithm~\ref{alg:smp} and Theorem~\ref{th:smp}.
%
\begin{small}
\begin{algorithm}
\KwData{$A \in \expander$; $y \in \R^m$}
 \KwResult{ $\hat x \in \R^n$ $\st$ $\|x - \hat x\|_1 = \mathcal{O}(\|y - Ax\|_1/d)$ }
 $\hat x \leftarrow 0$, $r \leftarrow y$\;
 \While{not converged}{
		$\hat x \leftarrow H_k\left[\hat x + H_{2k}[\mathcal{M}(r)]\right]$\;
		$r\leftarrow y - A\hat x$\;
}
\caption{SMP \citep{berinde2008practical}}
\label{alg:smp}
\end{algorithm}
\end{small}
%
\begin{theorem}[SMP \citep{berinde2008practical}]
\label{th:smp}
Let $A \in \expander$ 
 and let $y = Ax + \eta$ for $x \in \Sparse$. Then, there exists an $\varepsilon \ll 1$ such that SMP recovers $\hat x \in \R^n$ such that $\|x - \hat x\|_1 = \mathcal{O}(\|\eta\|_1/d)$. 
The algorithm terminates in $\mathcal{O}(\log(d\|x\|_1/\|\eta\|_1))$ iterations with complexity $\mathcal{O}(nd + n\log n)$.
\end{theorem}
%

%
\subsection{Sequential Sparse Matching Pursuit (SSMP)}

It was observed in \citep{berinde2009sequential} that SMP typically
failed to converge to the sought sparsest solution when the problem
parameters fall outside the region of theoretical guarantees.  
Though SMP updates each entry in $x$ to individually reduce the
$\ell_1$ norm of the residual, by updating multiple values of $x$ in
parallel causes SMP to diverge
even for
moderately small ratios of $k/m$.
%
To overcome these limitations, the authors proposed SSMP, which updates $\hat x$ sequentially rather than in parallel.
That is, at each iteration, SSMP will look for a single node $j \in [n]$ and an update $\omega \in \R$ that minimise $\|r - \omega a_j\|_1$, which can be found by computing $\argmax_{j \in [n]} \mathcal{M}(r)$, see the discussion in Section \ref{sec:median-minimise}.
This approach results in a strict decrease in $\|r\|_1$, but the
sequential update results in an overall increase in computational
complexity, see Table~\ref{table:ccs}.
SSMP and its theoretical guarantees are stated in
Algorithm~\ref{alg:ssmp} and Theorem~\ref{th:ssmp}.
%
\begin{small}
\begin{algorithm}
\KwData{$A \in \expander$; $y \in \R^m$; $c>1$;}
 \KwResult{ $\hat x \in \R^n$ $\st$ $\|x - \hat x\|_1 = \mathcal{O}(\|y - Ax\|_1/d)$ }
 $\hat x \leftarrow 0$, $r \leftarrow y$\;
 \While{not converged}{
			Find $(j, \omega) \in [n]\times \R$ $\st$ $\|r - \omega a_j\|_1$ is minimized\;
			$\hat x_j \leftarrow \hat x_j + \omega$\;
	        Perform $\hat x \leftarrow H_k[\hat x]$ every $(c-1)k$ iterations\;
			$r\leftarrow y - A\hat x$\;
}
\caption{SSMP \citep{berinde2009sequential}}
\label{alg:ssmp}
\end{algorithm}
\end{small}
\begin{theorem}[SSMP \citep{berinde2009sequential}]
\label{th:ssmp}
Let $A \in \mathbb{E}_{(c + 1)k, \varepsilon, d}$ and let $y = Ax +
\eta$ for $x \in \Sparse$. Then, there exists an $\varepsilon \ll 1$
such that SSMP with fixed $c>1$ recovers $\hat x \in \R^n$ such that $\|x - \hat x\|_1 = \mathcal{O}(\|\eta\|_1)$. The algorithm terminates in $\mathcal{O}(k)$ iterations of complexity $\mathcal{O}\left(\frac{d^3n}{m} + n + \left(\frac{n}{k}\log n\right)\log \|x\|_1\right)$.
\end{theorem}

\subsection{Left Degree Dependent Signal Recovery (LDDSR)}

LDDSR was proposed in \citep{xu2007efficient} and decodes by exploiting the information preserving property given in Lemma \ref{th:equivalence-lemma}.
The main insight is that one can lower bound the number of elements in $\{i \in [m]: y_i \in \argsupp(x)\}$, and use the structure of $A$ to find a $j \in [n]$ and a nonzero value $\omega \in \R$ that appears more than $d/2$ times in $r_{\mathcal{N}(j)}\in\RR^d$.
It is shown in \citep{xu2007efficient} that updating $\hat x_j \leftarrow \hat x_j + \omega$ guarantees a decrease in $\|r\|_0$ when $\varepsilon = 1/4$.
LDDSR and its theoretical guarantees are stated in
Algorithm~\ref{alg:lddsr} and Theorem~\ref{th:lddsr}.
%
\begin{small}
\begin{algorithm}
 \KwData{$A \in \expander$; $y \in \R^m$}
 \KwResult{ $\hat x \in \R^n$ $\st$ $ y = A \hat x$ }
 $\hat x \leftarrow 0$, $r \leftarrow y$\;
 \While{not converged}{
	Find $(j, \omega) \in [n]\times \R\setminus\{0\}$ $\st$ $|\{i \in \mathcal{N}(j): r_i = \omega\}| > \frac{d}{2}$\;
	$\hat x_j \leftarrow \hat x_j + \omega$\;
	$r \leftarrow y - A\hat x$\;
}
\caption{LDDSR \citep{xu2007efficient}}
\label{alg:lddsr}
\end{algorithm}
\end{small}
\begin{theorem}[LDDSR \citep{xu2007efficient}]
\label{th:lddsr}
Let $A \in \expander$ with $\varepsilon = 1/4$ and $x \in \Sparse$. Given $y = Ax$, LDDSR recovers $x$ in at most $\mathcal{O}(dk)$ iterations with complexity $\mathcal{O}(\frac{d^3n}{m} + n)$.
\end{theorem}
\subsection{Expander Recovery (ER)}
\label{sec:er}

ER \citep{jafarpour2008efficient} differs from LDDSR by considering
$\varepsilon\leq1/4$ and suitably adapting the set of indices from which
an entry in $\hat x$ may be updated.  This modification allows the
number of iterations guaranteed to be improved, see Theorem
\ref{th:lddsr}.  In particular, ER gurantees convergence in
$\mathcal{O}(k)$ iterations of complexity $\mathcal{O}(nd)$.
ER and its theoretical guarantees are stated in
Algorithm~\ref{alg:er} and Theorem~\ref{th:er}.
%
\begin{small}
\begin{algorithm}
 \KwData{$A \in \expander$; $y \in \R^m$}
 \KwResult{ $\hat x \in \R^n$ $\st$ $ y = A \hat x$ }
 $\hat x \leftarrow 0$, $r \leftarrow y$\;
 \While{not converged}{
	Find $(j, \omega) \in [n] \times \R\setminus\{0\}$ $\st$ $|\{i \in \mathcal{N}(j): r_i = \omega\}| \geq (1 - 2\varepsilon)d$.\;
	$\hat x_j \leftarrow \hat x_j + \omega$\;
	$r \leftarrow y - A\hat x$\;
}
\caption{ER \citep{jafarpour2008efficient}}
\label{alg:er}
\end{algorithm}
\end{small}
\begin{theorem}[ER \citep{jafarpour2008efficient}]
\label{th:er}
Let $A \in \mathbb{E}_{2k, \varepsilon, d}$ with $\varepsilon \leq 1/4$ and $m = \mathcal{O}(k \log (n/k))$. Then, for any $x \in \Sparse$, given $y = Ax$, ER recovers $x$ in at most $\mathcal{O}(k)$ iterations of complexity $\mathcal{O}(\frac{d^3 n}{m} + n)$.
\end{theorem}
Though ER seemingly requires knowledge of $\varepsilon$ to implement, which is
NP-hard to compute, knowledge of $\varepsilon$ 
can be circumvented by selecting the node to update by 
\begin{equation}
\argmax_{j \in [n]} \left|\{ i \in \mathcal{N}(j) : r_i = \mode( r_{\mathcal{N}(j)})\}\right|.
\end{equation}
\subsection{Discussion}
\label{sec:discussion_prior_art}
Having introduced these algorithms, we now point out some important commonalities between them.
\subsubsection{Iterative greedy algorithms}
The CCS algorithms we have presented share the structure of Algorithm~\ref{alg:ccs}.

\begin{small}
\begin{algorithm}
 \KwData{$A \in \R^{m \times n}$; $y \in \R^m$}
 \KwResult{ $\hat x \in \R^n$ $\st$ $ y = A \hat x$ }
 $\hat x \leftarrow 0$, $r \leftarrow y$\;
  \While{not converged}{
  	Compute a {\em score} $s_j$ and an {\em update} $u_j$ $\forall$ $j \in [n]$\; 
	Select $S \subset [n]$ based on a rule on $s_j$\;
	$\hat x_j \leftarrow \hat x_j + u_j$ for $j \in S$\;
	$k$-threshold $\hat x$\;
	$r\leftarrow y - A\hat x$\;
}
\caption{Iterative greedy CCS algorithms}
\label{alg:ccs}
\end{algorithm}
\end{small}
%
%
The dominant computational cost in CCS greedy algorithms is
concentrated in computing $s_j$ and $u_j$, and in selecting the set
$S$ of nodes that will be updated. At each step of these algorithms, a
subset $S \subset [n]$ is selected. In SMP, we have $S = [n]$ which
makes it of sublinear complexity in $\|x\|_1$, but typically diverges
for even moderate values of $\rho:=k/m$. All other algorithms update a
single entry of $\hat x$ per iteration; that is, they choose $S
\subset [n]$ with $|S| = 1$. This brings benefits in terms of
convergence and recovery region, but compromises the computational
complexity of the algorithms.  A summary of these properties is given
in Table~\ref{table:ccs}.
%
%
\subsubsection{Median minimises $\|r\|_1$}
\label{sec:median-minimise}
The operation $\median( r_{\mathcal{N}(j)})$ can be recast as the problem of finding the scalar $\omega \in \R$ that minimises $\|r - \omega a_j\|_1$. To see this, note that the function
\begin{equation}
\| r - \omega a_j\|_1 = \sum_{i \in \mathcal{N}(j)} |r_{i} - \omega| + \mbox{constant}
\end{equation}
is at a minimum when $|\{i \in \mathcal{N}(j) : r_i -\omega > 0\}| = |\{i \in \mathcal{N}(j) : r_i - \omega < 0\}|$. Then, by definition of the median,
\begin{equation}
\label{eq:median-minimiser}
\argmin_{\omega \in \R} \|r - \omega a_j\|_1 = \median( r_{\mathcal{N}(j)})
\end{equation}
This is independent of the expansion parameter $\varepsilon$.
\subsubsection{Mode does not minimise $\|r\|_0$}
In \citep{jafarpour2008efficient, xu2007efficient}, it is shown that Algorithms \ref{alg:lddsr} and \ref{alg:er} use Lemma \ref{th:equivalence-lemma} to find a pair $(j, \omega)$ such that
\begin{equation}
\label{eq:contraction_1}
\|y - A(\hat x + \omega e_j)\|_0 < \|y - A\hat x\|_0 - (1 - 4 \varepsilon)d.
\end{equation}
However, when $(y - A\hat x)_{\mathcal{N}(j)}$ does not contain any zeros, we can guarantee that  
\begin{equation}
\label{eq:contraction_2}
\|y - A(\hat x + \omega e_j)\|_0 < \|y - A\hat x\|_0 - (1 - 2 \varepsilon)d.
\end{equation}
For dissociated signals, where $\sum_{j \in \supp(x)} x_j \neq 0$, we
can always ensure that the greater contraction rate will be achieved. 
%
%
\subsubsection{Updating $s_j$ and $u_j$}

Algorithms \ref{alg:ssmp}, \ref{alg:lddsr}, \ref{alg:er} need to compute a score $s_j = s_j(r_{\mathcal{N}(j)})$ for each $j \in [n]$, which can be done at cost $\mathcal{O}(dn)$.
It is important to note that they do not need to recompute all the scores at each iteration. 
A common strategy is to compute each of the scores once and store them
with their corresponding node $j \in [n]$ in some data structure (like
priority queues \cite{berinde2009sequential} or red-black trees \cite{jafarpour2008efficient}).
Then, at each iteration, we can efficiently request the node $j \in [n]$ that maximises the score (median, mode, etc.) and use it to update $\hat x_j$.
This update will affect $d = |\mathcal{N}(j)|$ entries of the residual, so we only need to recompute the scores corresponding to $|\bigcup_{i \in \mathcal{N}(j)}\mathcal{N}(i)| = \mathcal{O}(d^2n/m)$ right nodes.

\begin{table*}[!htbp]
\begin{center}
\begin{tabular}{cc|c|c|c|c|c|c|}
\cline{3-8}
& & Objective & Score & Signal & Concurrency & Number of iterations & Iteration cost \\\cline{1-8}
\multicolumn{1}{ |c  }{\multirow{9}{*}{\rotatebox[origin=c]{90}{Prior art$\esp\esp$}}} &
\multicolumn{1}{ |c| }{} &   &  & &  &  &  \\
\multicolumn{1}{ |c  }{}                        &
\multicolumn{1}{ |c| }{SMP \citep{berinde2008practical}} & $\ell_1$  & median & any & parallel  & $\mathcal{O}(\log \|x\|_1)$ & $\mathcal{O}(nd + n\log n)$\\
\multicolumn{1}{ |c  }{}                        &
\multicolumn{1}{ |c| }{} &   & &   &  & &\\
\multicolumn{1}{ |c  }{}                        &
\multicolumn{1}{ |c| }{SSMP \citep{berinde2009sequential}} & $\ell_1$  & median & any & serial & $\mathcal{O}(k)$ & $\mathcal{O}(d^3n/m  + n +(\frac{n}{k}\log n)\log \|x\|_1)$\\
\multicolumn{1}{ |c  }{}                        &
\multicolumn{1}{ |c| }{} &   &  &   &  & &\\
\multicolumn{1}{ |c  }{}                        &
\multicolumn{1}{ |c| }{LDDSR \citep{xu2007efficient}} & $\ell_0$  & mode & any & serial & $\mathcal{O}(dk)$ & $\mathcal{O}(\frac{d^3n}{m} + n)$\\
\multicolumn{1}{ |c  }{}                        &
\multicolumn{1}{ |c| }{} &   &  &   &  & &\\
\multicolumn{1}{ |c  }{}                        &
\multicolumn{1}{ |c| }{parallel-LDDSR} & $\ell_0$  & mode & dissociated & parallel & $\mathcal{O}(\log k)$ & $\mathcal{O}(nd)$\\
\multicolumn{1}{ |c  }{}                        &
\multicolumn{1}{ |c| }{} &   &  &   &  & &\\
\multicolumn{1}{ |c  }{}                        &
\multicolumn{1}{ |c| }{ER \citep{jafarpour2008efficient}}  & $\ell_0$ & mode  & any & serial & $\mathcal{O}(k)$ & $\mathcal{O}(\frac{d^3n}{m}  + n)$\\
\multicolumn{1}{ |c  }{}                        &
\multicolumn{1}{ |c| }{} &   & &   &  & &\\
\hline
\hline
\multicolumn{1}{ |c  }{\multirow{7}{*}{\rotatebox[origin=c]{90}{$\;\;\;\esp$Contributions}}} &
\multicolumn{1}{ |c| }{} &   &   &  &  & &\\
\multicolumn{1}{ |c  }{}                        &
\multicolumn{1}{ |c| }{serial-$\ell_0$} & $\ell_0$ & $\ell_0$-decrease
& dissociated & serial & $\mathcal{O}(n\log k)$ & $\mathcal{O}(d)$\\
\multicolumn{1}{ |c  }{}                        &
\multicolumn{1}{ |c| }{} &   & &   &  & &\\
\multicolumn{1}{ |c  }{}                        &
\multicolumn{1}{ |c| }{parallel-$\ell_0$} & $\ell_0$  & $\ell_0$-decrease & dissociated & parallel & $\mathcal{O}(\log k)$ & $\mathcal{O}(nd)$ \\
\multicolumn{1}{ |c  }{}                        &
\multicolumn{1}{ |c| }{} &   &   &  &  & &\\
\hline
\end{tabular}
\end{center}
\caption{Summary of prior art in combinatorial compressed-sensing.}
\label{table:ccs}
\end{table*}

\section{Main contributions: Iterative $\ell_0$-minimisation}
\label{sec:contributions}

Our main contributions, Serial-$\ell_0$ and Parallel-$\ell_0$, 
advance combinatorial compressed sensing by
    having comparatively high phase transitions while retaining the low
    computational complexity of SMP and the parallel implementation of
    LDDSR.  In particular, Parallel-$\ell_0$
    is observed to typically recover the
      sparsest solution of underdetermined systems of equations in less
      time than any other compressed sensing algorithm when the signal is
      dissociated and the sensing matrix is an expander graph.  

Serial-$\ell_0$ and Parallel-$\ell_0$
look for a solution
 by identifying nodes which if updated sequentially would 
strictly reduce the $\|r\|_0$ by at least $\alpha$. That is, they will choose
a coordinate $j$ of $x$, and an update value $\omega$ such that,  
\begin{equation}
\label{eq:selection_l0}
{(j, \omega) \in [n]\times \R} \;\;\st\;\; \|r\|_0 - \|r - \omega a_j\|_0 \geq \alpha,
\end{equation}
for some $\alpha \in (1, d]$. By selecting a pair $(j, \omega)$ satisfying (\ref{eq:selection_l0}), 
Serial-$\ell_0$ yields a decrease in $\|r\|_0$ at every update, and is guaranteed to converge
in $\mathcal{O}(n\log k)$ iterations of computational complexity
$\mathcal{O}(d)$ if the signal is dissociated.  Parallel-$\ell_0$ is
designed similarly, but adapted to be able to take full
advantage of modern massively parallel computational resources.
Indeed, Parallel-$\ell_0$ selects and update all pairs $(j, \omega)$
satisfying (\ref{eq:selection_l0}) and updates these values in $x$ in parallel.
Under this updating scheme, a strict contraction in
$\|r\|_0$ is guaranteed at every iteration when the signal is dissociated and $\alpha = (1 -
2\varepsilon)d$ with $\varepsilon\leq1/4$, though we show in 
Section \ref{sec:numerical-experiments} that one can fix $\alpha = 2$ and get high phase transitions and
exceptional speed.

Section \ref{subsec:lemmas} presents the key technical lemmas that explain
the behaviour of an iteration of Serial-$\ell_0$ and
Parallel-$\ell_0$.  In particular, technical lemmas are stated to show
how often values in $Ax$ appear when $x \in \Sparse$ and
$A\in\expander$, and that when a value in $Ax$ appears sufficiently
often it must be a value from $x$ at a specified location.  This
property ensures the algorithm updates its approximation $\hat x$ with
values $x_j$ in the $j^{th}$ entry, that is with the exact
values from $x$ at the correct locations.  The dissociated signal
model, Definition \ref{def:dissociated-signals}, is an essential
component in the analysis presented in Section \ref{subsec:lemmas},
though we will observe that
the algorithms' recovery region degrade gracefully as the fraction of duplicate entries
in $x$ increases.  The convergence rate of Serial-$\ell_0$ and
Parallel-$\ell_0$ are presented in Section \ref{subsec:rate}, and
together they establish Theorem~\ref{th:convergence-expl0de}.

\subsection{Technical lemmas}\label{subsec:lemmas}
%
\begin{lemma}[Properties of dissociated signals]
\label{th:dissociated-signals}
Let $x \in \Sparse$ be dissociated. Then, 
\begin{enumerate}[label=(\roman*)]
\item $x_i \neq x_j$ $\;\forall\; i,j \in \supp(x), \;\; i\neq j$.
\item $\sum_{j \in T} x_j \neq 0$ $\;\forall\; \emptyset \neq T \subset \supp(x)$.
\end{enumerate}
\end{lemma}
\begin{IEEEproof}
The result follows from (\ref{eq:dissociated}). For (i) we set $T_1 = \{i\}$ and $T_2 = \{j\}$, and for (ii) we let $T_2 = \emptyset$.
\end{IEEEproof}

\begin{lemma}[Bounded frequency of values in expander measurements of dissociated signals]
\label{th:bounded_frequency}
Let $x \in \Sparse$ be dissociated, $A\in\expander$, and $\omega$ a nonzero value
in $Ax$. Then, there is a unique set $T \subset \supp(x)$ such that $\omega=\sum_{j\in T} x_j$ and the value $\omega$ occurs in $y$ at most $d$ times, 
\begin{equation}
\left|\left\{ i \in [m] : y_i = \omega\right\}\right| \leq d \;\;\;\forall \;\omega \ne 0.
\end{equation}
\end{lemma}
\begin{IEEEproof}
The uniqueness of the set $T \subset \supp(x)$ such that $\omega=\sum_{j\in T} x_j$ follows by the definition of dissociated. 
Since $|\mathcal{N}(j)| = d$ for all $j \in [n]$, we have that,
\begin{equation}
\left|\left\{ i \in [m] : y_i = \omega\right\}\right| = \left|\bigcap_{j\in T} \mathcal{N}(j)\right| \leq |\mathcal{N}(j_0)| = d
\end{equation}
for any $j_0\in T$.
\end{IEEEproof}

%
\begin{lemma}[Pairwise column overlap]
\label{lemma:2-intersections}
Let $A \in \expander$. If $\varepsilon \leq 1/4$, every pair of columns
of $A$ intersect in less than $(1 - 2\varepsilon)d$ rows, that is, for all $j_1, j_2 \in [n]$ with $j_1 \neq j_2$
\begin{equation}
\label{eq:two-sets}
\left|\mathcal{N}(j_1)\bigcap\mathcal{N}(j_2)\right|<(1-2\varepsilon)d.
\end{equation}
\end{lemma}
\begin{IEEEproof}
Let $S \subset [n]$ be such that $|S| = 2$ then 
\begin{equation}
\label{eq:2-intersections}
|\mathcal{N}(S)| > 2(1 - \varepsilon)d \geq 2d - (1 - 2\varepsilon)d,
\end{equation}
where the first inequality is Definition \ref{def:expander} and the
second inequality follows from $\varepsilon\leq1/4$.
However, $|\mathcal{N}(S)|$ can be rewritten as
\begin{equation}
\label{eq:NS-split}
\left|\mathcal{N}(S)\right| = \left|\mathcal{N}(j_1)\right| + \left|\mathcal{N}(j_2)\right| - \left|\mathcal{N}(j_1) \bigcap \mathcal{N}(j_2)\right|,
\end{equation}
for some $j_1, j_2 \in [n]$. Coupling \eqref{eq:NS-split} with \eqref{eq:2-intersections} gives \eqref{eq:two-sets}.
\end{IEEEproof}

\begin{lemma}[Progress]
\label{th:progress}
Let $y = Ax$ for dissociated $x \in \Sparse$ and $A \in \expander$ with $\varepsilon \leq 1/4$. There is a pair $(j, \omega) \in [n]\times \R$ such that 
\begin{equation}
\label{eq:identification}
|\{i \in \mathcal{N}(j) : y_i = \omega\}| \geq (1 - 2\varepsilon)d.
\end{equation}
\end{lemma}
\begin{IEEEproof}
Let $S = \supp(x)$, then by the information-preserving property
\eqref{eq:unique-neighbors} it holds that $|\mathcal{N}_1(S)| > (1 -
2\varepsilon)d|S|$, where $\mathcal{N}_1(S)$ is defined in
\eqref{eq:N1}, or alternatively, by $\mathcal{N}_1(S) = \{i \in [m] :
y_i = x_j, j \in S\}$ in the context of dissociated signals. Given the
lower bound in $|\mathcal{N}_1(S)|>(1 - 2\varepsilon)d|S|$, if $|S|\ne
0$, at least one $j \in S$ must have at least $(1 - 2\varepsilon)d$ neighbours in $y$ with identical nonzero entries. Letting $\omega$ take the value of such repeated nonzeros in $y$ gives the required pair $(j, \omega) \in [n]\times \R$.
\end{IEEEproof}

\begin{lemma}[Support identification]
\label{th:supp_identification}
Let $y = Ax$ for dissociated $x \in \Sparse$ and $A \in \expander$ with $\varepsilon \leq 1/4$. Let $\omega \neq 0$ be such that \eqref{eq:identification} and $\omega = x_j$.
\end{lemma}
\begin{IEEEproof}
Our claim is that for any $\omega$ which is a nonzero value from $y$,
if the cardinality condition \eqref{eq:identification} is satisfied
then the value $\omega=\sum_{j\in T}x_j$ occurs for the set $T$ being
a singleton, $|T|=1$.  Lemma \ref{th:bounded_frequency} states 
that $T$ is unique and that
\begin{equation}
|\{i \in \mathcal{N}(j) : y_i = \omega\}| =\left| \bigcap_{j\in T}\mathcal{N}(j) \right|.
\end{equation}
If $|T|>1$ then the above is not more than the cardinality of the intersection of any two of the sets
$\mathcal{N}(j_1)$ and $\mathcal{N}(j_2)$, and by \eqref{eq:two-sets} in Lemma
\ref{lemma:2-intersections} that is less than $(1-2\varepsilon)d$ which
contradicts the cardinality condition \eqref{eq:identification} and
consequently $|T|\leq1$. However, Lemma \ref{th:progress} guarantees that $|T|>0$, so $|T| = 1$ and $\omega=x_j$.
\end{IEEEproof}

Equipped with Lemmas \ref{th:dissociated-signals}
- \ref{th:supp_identification} we prove Theorem
\ref{th:convergence-expl0de} considering Serial-$\ell_0$ and
Parallel-$\ell_0$ separately, beginning with the later.
Note that since $x \in \Sparse$ and the algorithm only sets entries in $\hat x$ to the
correct values of $x$, then $x - \hat x \in \Sparse$, and Lemmas \ref{th:progress} and \ref{th:supp_identification} hold with $y$ replaced by $r = y - A(x - \hat x)$.
\subsection{Proof of Theorem \ref{th:convergence-expl0de}}\label{subsec:rate}
%
%
\begin{theorem}[Convergence of Parallel-$\ell_0$]
\label{th:convergence-parallel-l0}
Let $A \in  \expander$ and let $\varepsilon \leq 1/4$, and $x \in \Sparse$ be dissociated. Then, Parallel-$\ell_0$ with $\alpha = (1 - 2\varepsilon)d$ can recover $x$ from $y = Ax \in \R^m$ in $\mathcal{O}(\log k)$ iterations of complexity $\mathcal{O}(dn)$.
\end{theorem}
\begin{IEEEproof}
Let $\hat x = 0$ be our initial approximation to $x \in \Sparse$.
During the $\ell^{th}$ iteration of Parallel-$\ell_0$,
let $S_{\ell} = \supp(x-\hat{x})$ and include a subscript on the identification set $T=T_{\ell} \subset [n]$. 
As $A\in \expander$ and $\varepsilon \leq 1/4$,  
by Lemma \ref{th:supp_identification} and the required entry-wise
reduction in the residual by at least $\alpha=(1-2\varepsilon)d$, it
follows that Parallel-$\ell_0$ only sets entries in $\hat x$ to the
correct values of $x$ and as a result $\|x-\hat{x} \|_0\le \|x\|_0=k$
for every iteration. Moreover, by Lemma \ref{th:progress}, the set $T_\ell \neq \emptyset$ as long as $x \neq \hat x$, so the algorithm eventually converges.

In fact, we show that the rate of reduction of
$\|x-\hat{x} \|_0$ per iteration is by at least a fixed fraction
$\frac{2\varepsilon d}{1+\lfloor 2\varepsilon d \rfloor}$.  
As $A\in\expander$ has $d$ nonzeros per column, the
reduction in the cardinality 
of the residual, say $\|r^{\ell}\|_0-\|r^{\ell+1}\|_0$, can be at most
$d|T_{\ell}|$. That is,
\begin{equation}
\label{eq:bound-1}
\|r^{\ell}\|_0-\|r^{\ell+1}\|_0 \leq d|T_{\ell}|.
\end{equation}
To establish a fractional decrease in $|S_{\ell+1}|$ we develop a lower bound on $\|r^{\ell}\|_0-\|r^{\ell+1}\|_0$.
For $Q \subset S_\ell$ define the set $\mathcal{N}_1^{S_\ell}(Q)$ to be the set of nodes in $\mathcal{N}_1(S_\ell)$ and such that $i \in \mathcal{N}(j)$ for some $j \in Q$, $i.e.$ 
\begin{equation}
\label{eq:N1-restricted}
\mathcal{N}_1^{S_\ell}(Q) = \{i \in \mathcal{N}_1(S_\ell) : i \in \mathcal{N}(j), j \in Q\}. 
\end{equation}
Consider the partition $S_\ell = T_\ell \cup (S_\ell \setminus T_\ell)$ and rewrite $\mathcal{N}_1(S_\ell)$ as the disjoint union
\begin{equation}
\mathcal{N}_1(S_\ell) = \mathcal{N}_1^{S_\ell}(T_\ell) \cup \mathcal{N}_1^{S_\ell}(S_\ell \setminus T_\ell).
\end{equation}
Note that $ \mathcal{N}_1^{S_\ell}(T_\ell) \neq \mathcal{N}_1(T_\ell)$, and that by \eqref{eq:N1-restricted} and the dissociated signal model, $\mathcal{N}_1^{S_\ell}(T_\ell) \subset [m]$ is the set of indices in $r^{\ell}$ that are identical to a nonzero in $x$ and that have a frequency of at least $\alpha = (1 - 2\varepsilon)d$, so
\begin{equation}
\label{eq:res-lower-N1}
\|r^{\ell}\|_0-\|r^{\ell+1}\|_0 \geq |\mathcal{N}_1^{S_\ell}(T_\ell)|.
\end{equation}
At iteration $\ell$, if $T_\ell = S_\ell$, the full support of $x$ is correctly identified, so $x = \hat x$ after updating $\hat x$.
Otherwise, $T_{\ell} \neq S_{\ell}$ and the set $S_\ell \setminus T_\ell$ is not identified by the algorithm at this iteration.
We derive a lower bound on $|\mathcal{N}_1^{S_\ell}(T_\ell)|$ by considering two cases: $\alpha \in \mathbb{N}$ and $\alpha \notin \mathbb{N}$.

If $\alpha \in \mathbb{N}$, then each node in $S_\ell \setminus T_\ell$ has at most $\alpha - 1$ duplicates in $r^{\ell}$, so
\begin{equation}
\label{eq:bound-N1-SmT}
\left| \mathcal{N}_1^{S_\ell}(S_{\ell}\setminus T_{\ell})\right| \leq \left( \alpha -1\right)|S_{\ell} \setminus T_{\ell}|.
\end{equation}
Using the the information-preserving property \eqref{eq:unique-neighbors} and the identity given in \eqref{eq:N1-restricted} it follows that
\begin{align}
& |\mathcal{N}_1^{S_\ell}(T_\ell)| + |\mathcal{N}_1^{S_\ell}(S_\ell \setminus T_\ell)| \nonumber \\ 
&>(1 - 2\varepsilon)d\left(\left|T_\ell| + |S_\ell \setminus T_\ell\right|\right)\nonumber \\ 
\label{eq:bound-part}
&=(1 - 2\varepsilon)d|T_\ell| + |S_\ell \setminus T_\ell| + (\alpha - 1)|S_\ell \setminus T_\ell|.
\end{align}
Now, using \eqref{eq:bound-N1-SmT} to lower bound \eqref{eq:bound-part}, and solving for $|\mathcal{N}_1^{S_\ell}(T_\ell)|$ gives
\begin{equation}
\label{eq:lower-bound-N1T}
|\mathcal{N}_1^{S_\ell}(T_\ell)| \geq (1 - 2\varepsilon)d|T_\ell| + |S_\ell \setminus T_\ell|.
\end{equation}
By coupling \eqref{eq:lower-bound-N1T}, \eqref{eq:res-lower-N1}, and \eqref{eq:bound-1} into a chain of inequalities it is seen that
\begin{equation}
(1 - 2\varepsilon)d |T_{\ell}|+(|S_{\ell}|-|T_{\ell}|)\leq d|T_{\ell}|,
\end{equation}
which simplifies to
\begin{equation}
\label{eq:TvsS}
|T_{\ell}|\ge \frac{1}{1+2\varepsilon d}|S_{\ell}|.
\end{equation}

If $\alpha \notin \mathbb{N}$, then each node in $S_\ell \setminus T_\ell$ has at most $\lfloor\alpha\rfloor$ duplicates in $r^{\ell}$, so
\begin{equation}
\label{eq:bound-N1-SmT-notin}
\left| \mathcal{N}_1^{S_\ell}(S_{\ell}\setminus T_{\ell})\right| \leq \lfloor \alpha \rfloor|S_{\ell} \setminus T_{\ell}|.
\end{equation}
Similarly as in the former case, using \eqref{eq:N1-restricted} and the the information-preserving property \eqref{eq:N1}, we obtain
\begin{align}
& |\mathcal{N}_1^{S_\ell}(T_\ell)| + |\mathcal{N}_1^{S_\ell}(S_\ell \setminus T_\ell)| \nonumber \\ 
&>(1 - 2\varepsilon)d\left(\left|T_\ell| + |S_\ell \setminus T_\ell\right|\right)\nonumber \\ 
\label{eq:bound-part-notin}
&=(1 - 2\varepsilon)d|T_\ell| + (\alpha - \lfloor \alpha \rfloor)|S_\ell \setminus T_\ell| + \lfloor\alpha\rfloor|S_\ell \setminus T_\ell|.
\end{align}
Just as in the previous case, \eqref{eq:bound-part-notin} is bounded from below using \eqref{eq:bound-N1-SmT-notin}, and the resulting inequality is used to get
\begin{equation}
\label{eq:lower-bound-N1T-notin}
|\mathcal{N}_1^{S_\ell}(T_\ell)| \geq (1 - 2\varepsilon)d|T_\ell| + (\alpha - \lfloor \alpha \rfloor)|S_\ell \setminus T_\ell|.
\end{equation}
Inequalities \eqref{eq:lower-bound-N1T-notin}, \eqref{eq:res-lower-N1}, and \eqref{eq:bound-1} are then used to derive
\begin{equation}
\label{eq:chain-res-notin}
\alpha |T_{\ell}|+(\alpha - \lfloor \alpha \rfloor)(|S_{\ell}|-|T_{\ell}|)\leq d|T_{\ell}|.
\end{equation}
It follows from $\alpha = (1 - 2\varepsilon)d \notin \mathbb{N}$ and the properties of step functions that
\eqref{eq:chain-res-notin} is equivalent to
\begin{equation}
\label{eq:TvsS-notin-expanded}
|T_{\ell}|\ge \frac{1 - 2\varepsilon d + \lfloor 2 \varepsilon d \rfloor}{1 +  \lfloor 2 \varepsilon d \rfloor}|S_{\ell}|.
\end{equation}
Finally, note that \eqref{eq:TvsS-notin-expanded} reduces to \eqref{eq:TvsS} when $\alpha \in \mathbb{N}$, so using $S_{\ell+1}=S_{\ell}\setminus T_{\ell}$ and \eqref{eq:TvsS-notin-expanded}, we conclude that 
\begin{equation}
\label{eq:contraction-notint}
|S_{\ell+1}|\le \frac{2 \varepsilon d}{1 + \lfloor 2 \varepsilon d\rfloor}|S_{\ell}|.
\end{equation}

Since $|S_0|=k$ it follows that Parallel-$\ell_0$ will have converged after $\ell^*$
iterations when  
$k({2\varepsilon d}/({1 + \lfloor 2\varepsilon d\rfloor}))^{\ell^*}<1$, which
is achieved for
\begin{equation}
\ell^*\geq \left(\log^{-1} \left(\frac{1 + \lfloor 2\varepsilon d\rfloor}{2\varepsilon d}\right)\right)\log k.
\end{equation}
Each iteration of Parallel-$\ell_0$ involves computing
\eqref{eq:identification} for each $j\in [n]$, which is equivalent to
$n$ instances of finding the mode of a vector of length $d$ which can
be solved in $\mathcal{O}(d)$ complexity provided $\alpha>\lfloor d/2\rfloor$\cite{boyer1991mjrty}.
\end{IEEEproof}

\begin{theorem}[Convergence of Serial-$\ell_0$]
\label{th:convergence-sequential-l0}
Let $A \in  \expander$ and let $\varepsilon \leq 1/4$, and $x \in \Sparse$ be a dissociated signal. Then, Serial-$\ell_0$ with $\alpha = (1 - 2 \varepsilon)d$ can recover $x$ from $y = Ax \in \R^m$ in $\mathcal{O}(n\log k)$ iterations with complexity $\mathcal{O}(d)$.
\end{theorem}
\begin{IEEEproof}
The loop over $j\in [n]$ for Serial-$\ell_0$ identifies singletons $T$
to update values in $\hat{x}$ in serial.  The union of the singletons
for $j\in [n]$ includes the set of all nodes for which the
residual would be reduced by at least $\alpha$ if one were to forgo
the serial update in $\hat{x}$.  For $\alpha=(1-2\varepsilon)d$, the
proof of convergence for Theorem \ref{th:convergence-parallel-l0} establishes
that this results in a reduction of the cardinality of
$\supp(x-\hat{x})$ by at least a fraction
$2\varepsilon d/(1+\lfloor 2\varepsilon d\rfloor)$.  That is, for $p$ an integer,
Serial-$\ell_0$ satisfies
\begin{equation}
|\supp(x-\hat{x})|\le k\left(\frac{2\varepsilon d}{1+\lfloor 2\varepsilon d\rfloor}\right)^p
\end{equation}
after $\ell = p n$
iterations, and converges to $\hat{x}=x$ after at most $p^*>\log(k)/\log((1+\lfloor 2\varepsilon d\rfloor)/(2\varepsilon d))$ for convergence after
\begin{equation}
\ell^*\geq n \left(\log^{-1} \left(\frac{1 + \lfloor 2\varepsilon d \rfloor}{2\varepsilon d}\right)\right)\log k.
\end{equation}
iterations. Each iteration of Serial-$\ell_0$ involves computing the mode of a
vector of length $d$ and updating $d$ entries in the residual. Since we are interested in knowing the mode of $r_{\mathcal{N}(j)}$ only when the most frequent element occurs more than $d/2$ times,  this value can be found at cost $\mathcal{O}(d)$ \cite{boyer1991mjrty}.
\end{IEEEproof}

\subsection{Discussion}
\subsubsection{The computational cost of computing a mode can be improved if $d$ is small}
Evaluating (\ref{eq:identification}) for a given column $j \in
\supp(x)$ is equivalent to finding the mode of
$r_{\mathcal{N}(j)}$. This can be done at cost $\mathcal{O}(d)$ using
the Boyer-Moore Majority vote algorithm
\cite{boyer1991mjrty}. However, this algorithm requires that an
element of the array occurs more than $\lfloor d/2 \rfloor$ times, so
it might fail when we set $\alpha \in [\lfloor d/2 \rfloor]$. Our
numerical experiments (Section \ref{sec:numerical-experiments}) show
that best recovery regions are obtained for $\alpha = 2$, so we prefer
to have an algorithm with $\mathcal{O}(d)$ per-iteration cost for all
$\alpha \in [d]$.

Our approach is presented in Algorithm \ref{alg:score-expl0de}. Instead of looking for an $\omega \in \R$ satisfying (\ref{eq:identification}) for each $j \in [n]$, at the $\ell^{th}$ iteration we consider the reduction caused by $\omega_j$, defined as the $\ell \Mod{d}$-th element in $r_{\mathcal{N}(j)}$. When using this shifting strategy we compromise the final number of iterations, but we also keep a fixed cost of $d$ complexity per iteration for any $\alpha \in [d]$. The convergence guarantees of our algorithms when using this shifting strategy are presented in Theorem \ref{th:convergence-shifted-parallel-l0}.

\begin{small}
\begin{algorithm}
 \KwData{$j \in [n]$; $r \in \R^m$; $\omega \in \N$}
 \KwResult{$s_j \leftarrow |\{i \in \mathcal{N}(j): r_i = \omega\}|$}
\caption{Computation of score for serial-$\ell_0$ and parallel-$\ell_0$.}
\label{alg:score-expl0de}
\end{algorithm}
\end{small}

\begin{theorem}[Convergence of Shifted Parallel-$\ell_0$]
\label{th:convergence-shifted-parallel-l0}
Let $A \in  \expander$ with $\varepsilon \leq 1/4$, and $x \in \Sparse$ be dissociated. Then, the shifted versions of Serial-$\ell_0$ and Parallel-$\ell_0$ with $\alpha = (1 - 2\varepsilon)d$ can recover $x$ from $y = Ax \in \R^m$ in an average of $\mathcal{O}(dn\log k)$ operations.
\end{theorem}
\begin{IEEEproof}
Let $\hat x = 0$ be the initial approximation to $x \in \Sparse$, and
$A \in \expander$ with $\varepsilon \leq 1/4$. At $\ell^{th}$ iteration,
let $T = T_\ell$ be the set satisfying (\ref{eq:identification}), that
is, the one that Parallel-$\ell_0$ has marked for update. 
For $j \in T$, let $\omega_j$ be the most frequent element in $r_{\mathcal{N}(j)}$. In shifted-parallel-$\ell_0$, $\omega_j$ is not directly computed. Instead, at iteration $\ell$, the frequency of the $\ell \Mod{d}$-th value in $r_{\mathcal{N}(j)}$ is computed using Algorithm \ref{alg:score-expl0de} and tested against the imposed threshold $\alpha$. In the worst case, this increases the number of iterations by a factor $\mathcal{O}(d)$. However, on average, this is not the case, and convergence in $\mathcal{O}(\log k)$ iterations is guaranteed.

To see this, let $j \in T$ and let $\omega$ be drawn at random from $r_{\mathcal{N}(j)}$. Then, $\Pr(\omega = \omega_j) \ge 1 - 2\varepsilon$, so on average at iteration $\ell$ we will identify $|T_\ell|(1 - 2\varepsilon)$ correct entries in $\supp(x - \hat x)$. Given the bound for $|T_\ell|$ \eqref{eq:TvsS-notin-expanded} in the proof of parallel-$\ell_0$, we have that at each iteration we identify at least $(1 - 2\varepsilon)\frac{1 - 2\varepsilon d + \lfloor 2\varepsilon d\rfloor}{1 + \lfloor 2\varepsilon d \rfloor}|S_\ell|$. Therefore
\begin{equation}
|S_{\ell + 1}| \leq \left(\frac{(1 - 2\varepsilon)(1 - 2\varepsilon d + \lfloor 2 \varepsilon d\rfloor)}{ 1 + \lfloor 2 \varepsilon d \rfloor}\right) |S_{\ell}|.
\end{equation}
\end{IEEEproof}

\subsubsection{Our theoretical guarantees immediately apply to LDDSR}

When $\varepsilon = 1/4$, we have that $(1 - 2\varepsilon)d = d/2$, so we recover a parallel version of LDDSR (Algorithm \ref{alg:lddsr}) for dissociated signals. We call this algorithm Parallel-LDDSR, and we test its performance in Section \ref{sec:numerical-experiments}.

\subsubsection{Non-dissociated signals can be recovered with a dissociated $A$}
\label{sec:dissociated-A}

There are many signals models in which the dissociated condition does not hold. For instance, if $x$ is a binary signal or has integer-valued nonzeros. In this case, the sensing matrix $A$ can be modified to make the nonzero elements of $x$ identifiable by our algorithms. In particular, scaling each column of the matrix by $i.i.d.$ random numbers coming from a continuous distribution introduces enough information in $y$ for our algorithms to correctly identify $\supp(x)$.

\subsubsection{Expander matrices preserve information of dissociated signals}

We now discuss the concept of dissociated signals under an Information Theory viewpoint. To do this, suppose that $(X_1, \dots, X_k)$ is a vector of $k$ random variables associated with $\{x_1, \dots, x_k\} = \supp(x)$ and that $(X_1, \dots, X_k) \sim p$ for some distribution $p$ supported on a finite set. Note that condition (iii) in Definition \ref{def:dissociated-signals} implies that,
\begin{equation}
x_{i_1} + \cdots + x_{i_\ell} \neq x_{j_1} + \cdots + x_{j_\ell} \esp \mbox{for} \esp i_1 \neq j_1, \dots, i_\ell \neq j_\ell.
\end{equation}
Now, consider the following Shannon-entropy inequalities,
%
\begin{lemma}[Entropy inequalities]
\label{th:entropy-inequalities}
For a random variable $X \sim p$, let $H(\cdot)$ be its Shannon entropy. Now, let $X_1, \cdots, X_k$ be a set of random variables with joint distribution $(X_1, \dots, X_k) \sim p$. Assume that the random variable $X_i$ is supported on $\{(x_i)_1, \dots, (x_i)_\ell\}$. Then,
\begin{equation}
\label{eq:entropy-inequality}
H(X_1 + \cdots + X_k) \leq H(X_1, \dots, X_k) \leq H(X_1) + \cdots H(H_k)
\end{equation}
With equality on the left if and only if $(x_1)_{i_1} + \cdots + (x_k)_{i_k} \neq (x_1)_{j_1} + \cdots + (x_k)_{j_k}$ for $i_l \neq j_l$, and equality on the right if and only if $X_i \perp X_j$ for $i \neq j$.
\end{lemma}
\begin{IEEEproof}
See \citep{cover2012elements} and \citep{kelbert2013information} for a proof.
\end{IEEEproof}

In the case of discretely supported distributions, a dissociated signal can be understood as one in which the entries on $\supp(x)$ are drawn according to a distribution $p$ fulfilling,
\begin{enumerate}[label=(\roman*)]
\item $\Prob[X_i = \omega \mid X_j = \omega] = 0$ $\forall$ $i \neq j$ $\forall$ $\omega \neq 0$.
\item $\Prob[\sum_{j \in T} X_j = 0] = 0$ $\forall$ $T \subset [k]$
\item $\Prob[\sum_{j \in T_1} X_j = \sum_{j \in T_2} X_j] = 0$ $\forall$ $T_1, T_2 \subset [k]$ with $T_1 \neq T_2$.
\end{enumerate}
Property (iii) above, together with Lemma (\ref{th:entropy-inequalities}) say that probability distribution on the support of dissociated signals imply
\begin{equation}
\label{eq:information-preserving}
H\left(\sum_{j \in T} X_j\right) = H(X_1, \dots, X_k) \esp\forall\esp T\subset [k]
\end{equation}
And since the value of each entry in $y = Ax$ is distributed according to $\sum_{j \in T}X_j$ for some $T \subset [k]$, we get that when computing $y$ with a $A \in \expander$ having $\varepsilon \leq 1/4$ and a dissociated signal $x$, (\ref{eq:information-preserving}) will hold. This implies that  linear transformations with expander matrices preserve the information in $x$.
\section{Numerical Experiments}
\label{sec:numerical-experiments}
In this section we perform a series of numerical experiments to
compare Parallel-$\ell_0$ and Serial-$\ell_0$ with 
state-of-the-art compressed sensing algorithms.
These comparisons are done by adding Parallel-$\ell_0$ and
Serial-$\ell_0$ to the GAGA software package \cite{blanchard2012gpu}
which includes {\tt CUDA-C} implementations of a number of compressed
sensing algorithms as well as a testing environment to rapidly
generate synthetic problem instances.  This approach allows us to solve
hundreds of thousands of randomly generated problems and to solve
problems with $n$ in the millions.

Unless otherwise stated, all tests were performed with the nonzeros of
$x$ drawn from a standard normal distribution $\mathcal{N}(0,1)$ and
the parameter $\alpha$ in Serial-$\ell_0$ and Parallel-$\ell_0$ was
set to 2.

Figures \ref{fig:phase-transitions}-\ref{fig:iterations-linear} were computed using a Linux machine with Intel Xeon E5-2643 CPUs @ 3.30 GHz, NVIDIA Tesla K10 GPUs, and executed from Matlab R2015a.
Figures \ref{fig:timing_low_delta}-\ref{fig:parallel_l0_several_d} were computed using a Linux machine with Intel Xeon E5-2667 v2 CPUs @ 3.30GHz, NVIDIA Tesla K40 GPUs, and executed from Matlab R2015a.

\subsection{Substantially higher phase transitions}

\begin{figure}[!htbp]
	\centering
	\includegraphics[scale=0.39]{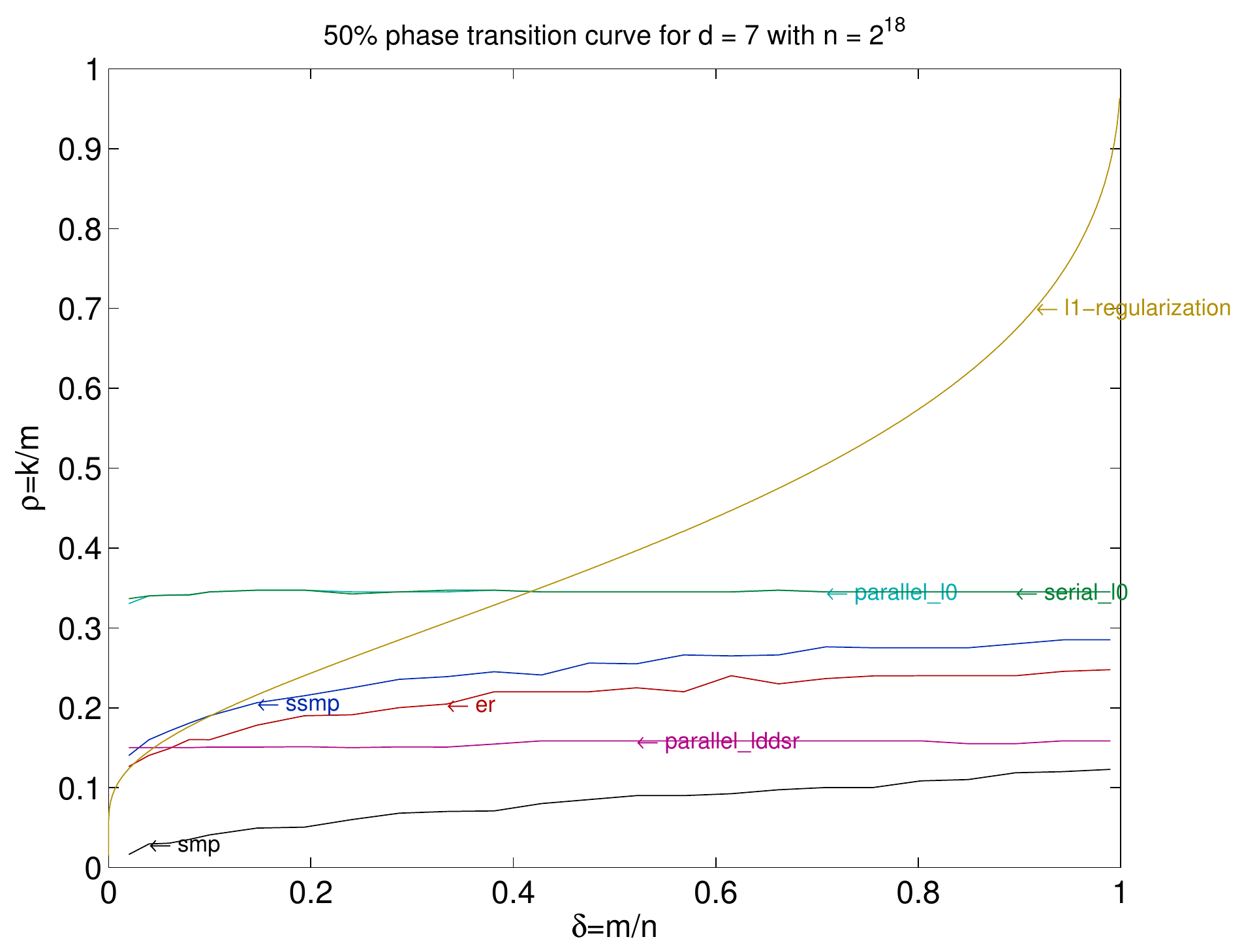}
	\caption[Phase Transitions for CCS algorithms]{50\% recovery probability logistic regression curves for $\mathbb{E}_{\varepsilon, 7, k}$ and $n = 2^{18}$. The curve for $\ell_1$-regularisation is the theoretical curve for dense Gaussian ensembles, and is shown for reference.}
	\label{fig:phase-transitions}
\end{figure}

The phase transition of a compressed-sensing algorithm
\citep{donoho2010precise} is the largest
value $k/m$, which we denote $\rho^*(m/n)$ noting its dependence on
$m/n$, for which the algorithm is typically (say greater than half of
the instances) able recovery all $k$ sparse vectors with
$k<m \rho^* (m/n)$.  The value $\rho^* (m/n)$ often converges to a fixed
value as $n$ is increased with $m/n$ being a fixed fraction.
Figure \ref{fig:phase-transitions} shows the phase transition curve for
each of the CCS algorithms stated in Section \ref{sec:prior-art}, as
well as Parallel-$\ell_0$ and Serial-$\ell_0$.  To facilitate
comparison with nonCCS algorithms, Figure \ref{fig:phase-transitions}
also includes the theoretical phase transition curve for
$\ell_1$-regularization for 
$A$ drawn Gaussian \cite{donoho2006high,donoho2005neighborly}, which
is observed to be consistent \cite{DTobserved} with
$\ell_1$-regularization for $A\in\expander$. 
The curves were computed by setting $n = 2^{18}$, $d= 7$, and a tolerance of $10^{-6}$.
The testing is done at $m=\delta_p n$ for 
\begin{equation*}
\delta_p \in \{ 0.02p :p \in [4]\} \cup \left\{0.1 + \frac{89}{1900}(p - 1) : p \in [20]\right\}.
\end{equation*}
For each $\delta_p$, we set $\rho = 0.01$ and generate 10 synthetic problems to be applied to the algorithms, with $x$ having independent and identically distributed normal Gaussian entries. With this restrictions, our signals are dissociated.
If at least one such problem was recovered successfully, we increase $\rho$ by $0.01$ and repeat the experiment.
The recovery data is then fitted using a logistic function in the
spirit of \citep{blanchard2013performance} and the 50\% recovery
transition of the logistic function is computed and shown in Figure~\ref{fig:phase-transitions}.

Note the low phase-transition curve of SMP and the substantially higher phase-transition curve of Parallel-$\ell_0$ and Serial-$\ell_0$.
As mentioned previously, the multiple updating mechanism of SMP gives it sublinear convergence guarantees, but greatly compromises its region of recovery.
We emphasise that the phase transition curves for Serial-$\ell_0$ and Parallel-$\ell_0$ are higher than those for SMP, SSMP, ER, and parallel-LDDSR. In particular, they are even higher than $\ell_1$-regularisation for $\delta \lessapprox 0.4$.

\subsection{Fastest compressed sensing algorithm}
\begin{figure}[!htbp]
	\centering
	\includegraphics[scale=0.39]{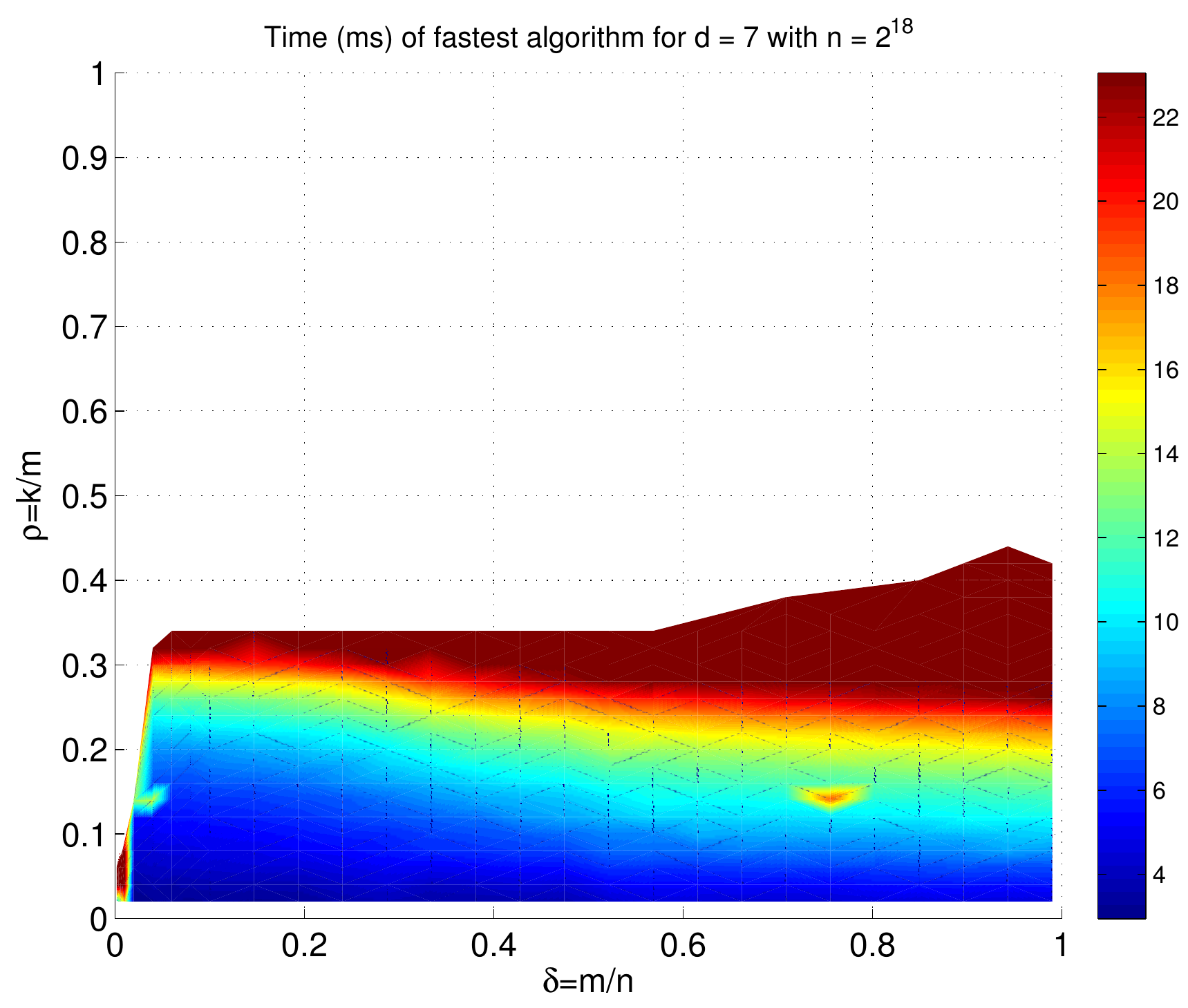}
	\caption[Time (ms) of fastest algorithm]{Average recovery time (ms) of the fastest algorithm at each $(\delta, \rho)$ for $\mathbb{E}_{k, \varepsilon, 7}$ and $n = 2^{18}$.}
	\label{fig:time-fastest}
\end{figure}

\begin{figure}[!htbp]
	\centering
	\includegraphics[scale=0.39]{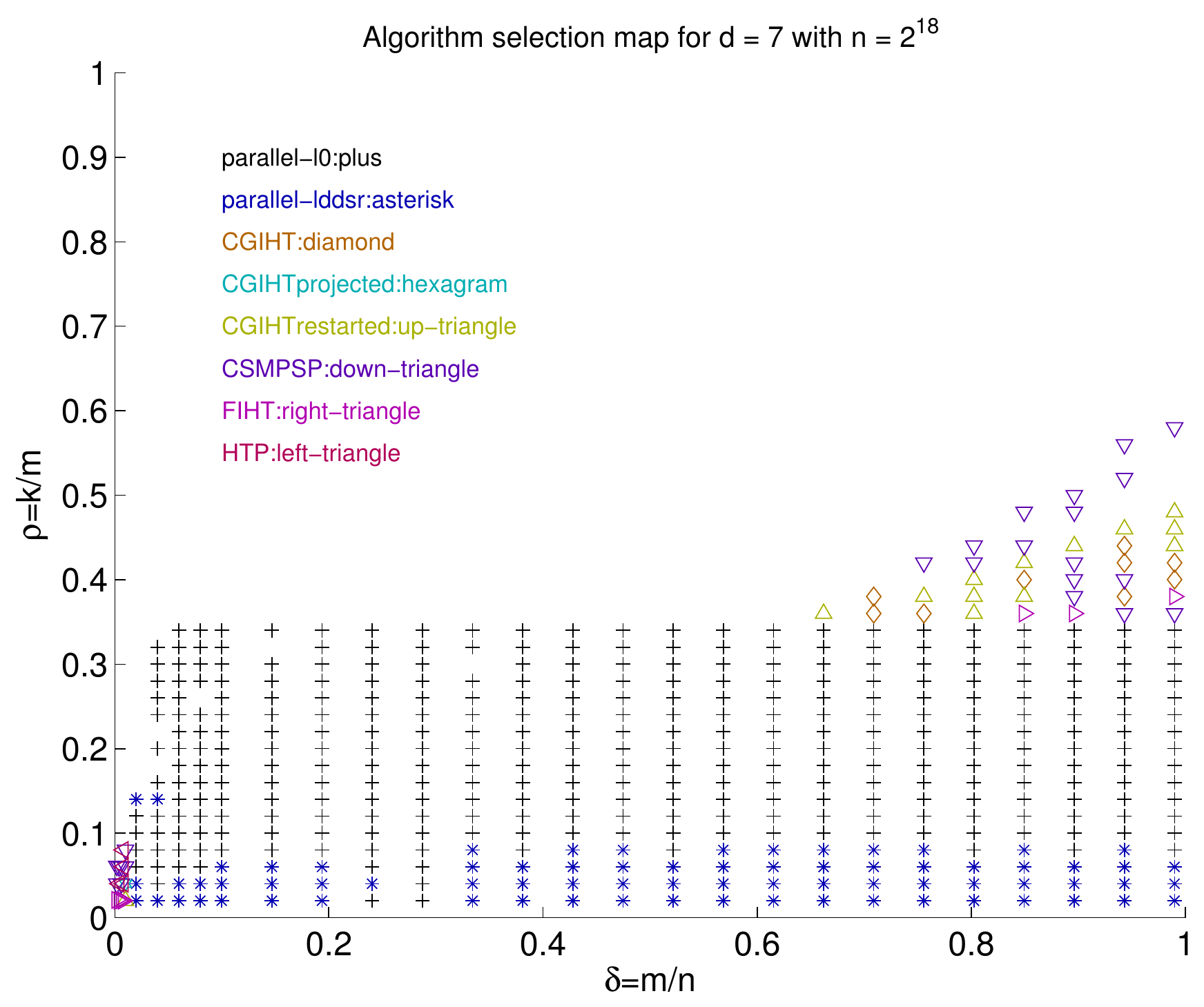}
	\caption[Algorithm selection map]{Selection map of the fastest algorithm at each $(\delta, \rho)$ for $\mathbb{E}_{k, \varepsilon, 7}$ and $n = 2^{18}$.}
	\label{fig:algorithm-selection-map}
\end{figure}

When the signal is dissociated, Parallel-$\ell_0$ is generally the
fastest algorithm for matrices $A\in\expander$.  We show this
numerically by computing the phase transitions of  
\begin{center}
{\small Serial-$\ell_0$, Parallel-$\ell_0$, parallel-LDDSR, ALPS, CGIHT, CSMPSP, ER, FIHT, HTP, NIHT, SMP, SSMP;}
\end{center}
and comparing their average time to convergence at each point of $(\delta, \rho)$.
The phase transitions are computed similarly to those in Figure \ref{fig:phase-transitions}, with problem parameters of $n = 2^{18}$ and $d = 7$. 
In particular, Parallel-$\ell_0$ is also used with $\alpha = 2$.
The results are shown in Figures \ref{fig:time-fastest} and \ref{fig:algorithm-selection-map}.
Specifically, Figure \ref{fig:time-fastest} shows the time in milliseconds that the fastest algorithm takes to converge when the problem parameters are located at $(\delta, \rho)$.
The fastest algorithm is in turn identified in Figure
\ref{fig:algorithm-selection-map}, where we can see that
Parallel-$\ell_0$ is consistently the fastest algorithm within its
phase transition, except for $\rho\ll 1$ where parallel-LDDSR takes
less time.
However, we note that the convergence guarantees of
parallel-LDDSR come as a byproduct of our analysis the domain in which
it is faster than Parallel-$\ell_0$ is the region of least importance
for applications as it indicates more than three fold more
measurements were taken than would have been necessary if
Parallel-$\ell_0$ were used.

\subsection{Parallelisation brings important speedups: examples with
  $m\ll n$}

\begin{figure}[!htbp]
	\centering
	\includegraphics[scale=0.39]{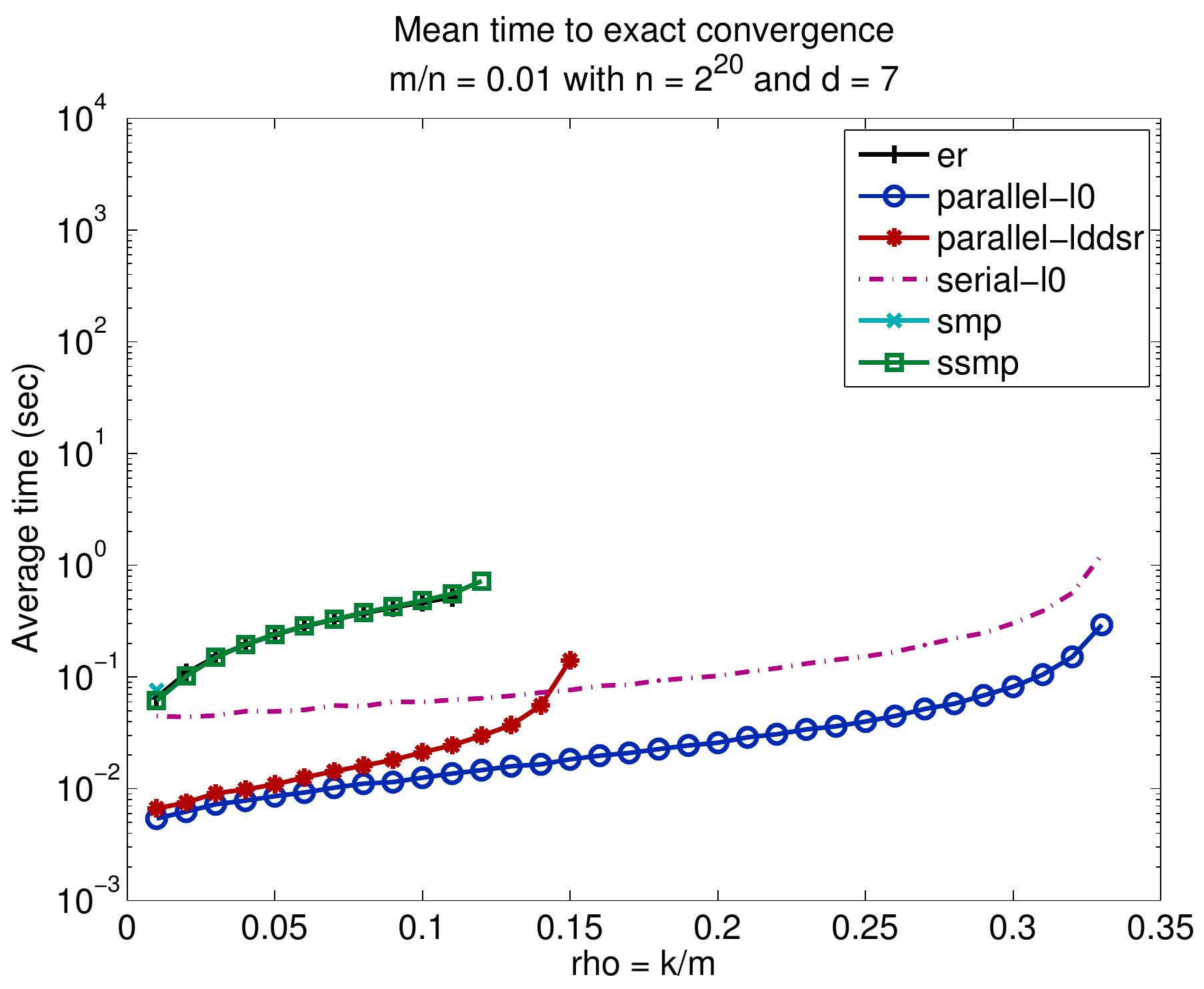}
	\caption[Mean time to exact convergence]{Average recovery time (sec) with dependence on $\rho$ for $\delta = 0.01$ and $\mathbb{E}_{k, \varepsilon, 7}$ with $n = 2^{20}$.}
	\label{fig:mean-time-exact-convergence-1}
\end{figure}

\begin{figure}[!htbp]
	\centering
	\includegraphics[scale=0.39]{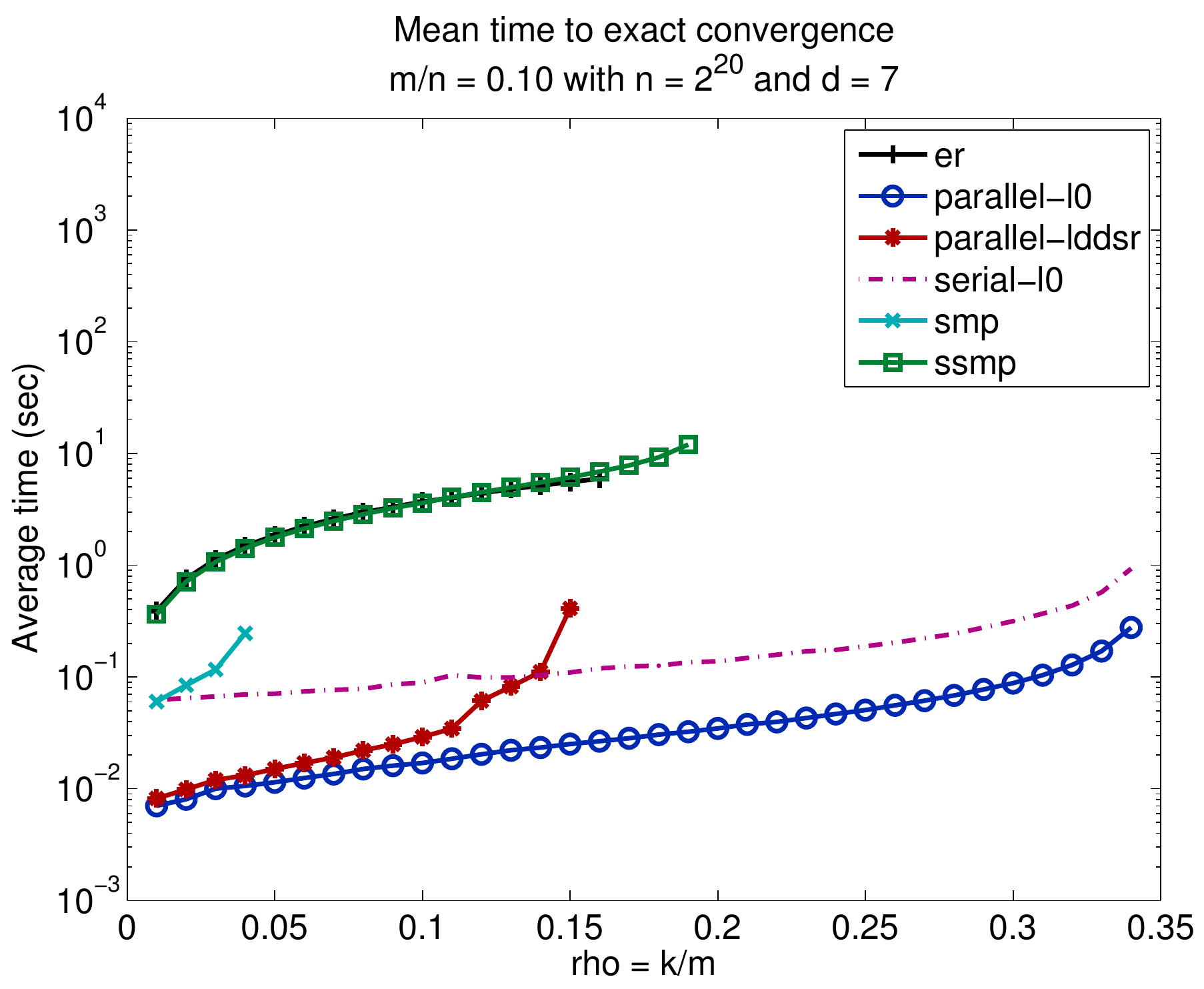}
	\caption[Mean time to exact convergence]{Average recovery time (sec) with dependence on $\rho$ for $\delta = 0.1$ and $\mathbb{E}_{k, \varepsilon, 7}$ with $n = 2^{20}$.}
	\label{fig:mean-time-exact-convergence-10}
\end{figure}

As shown in Algorithm \ref{alg:ccs}, the speed of Algorithms \ref{alg:smp}-\ref{alg:er} can be improved if the scores $s_j$ and updates $u_j$ are computed in parallel for each $j \in [n]$.
However, implementing this parallelisation is not enough to cut down an algorithm's complexity to that of the state-of-the-art's.
Figures \ref{fig:mean-time-exact-convergence-1}-\ref{fig:mean-time-exact-convergence-10} show the average time to exact convergence for each of the combinatorial compressed sensing algorithms.
It can be seen in addition to Serial-$\ell_0$ and Parallel-$\ell_0$
having higher phase transition than ER and SSMP, they are also
substantially faster to converge to the true solution for $n=2^{20}$
and either $\delta = 0.01$ or $\delta=0.1$.
It is interesting to note that for this problem size Serial-$\ell_0$
is substantially faster than ER and SSMP, even when the two latter are
implemented in parallel and run on a modern high performance computing
GPU. 

\subsection{Convergence in $\mathcal{O}(\log k)$ iterations}
\begin{figure}[!htbp]
	\centering
	\includegraphics[scale=0.39]{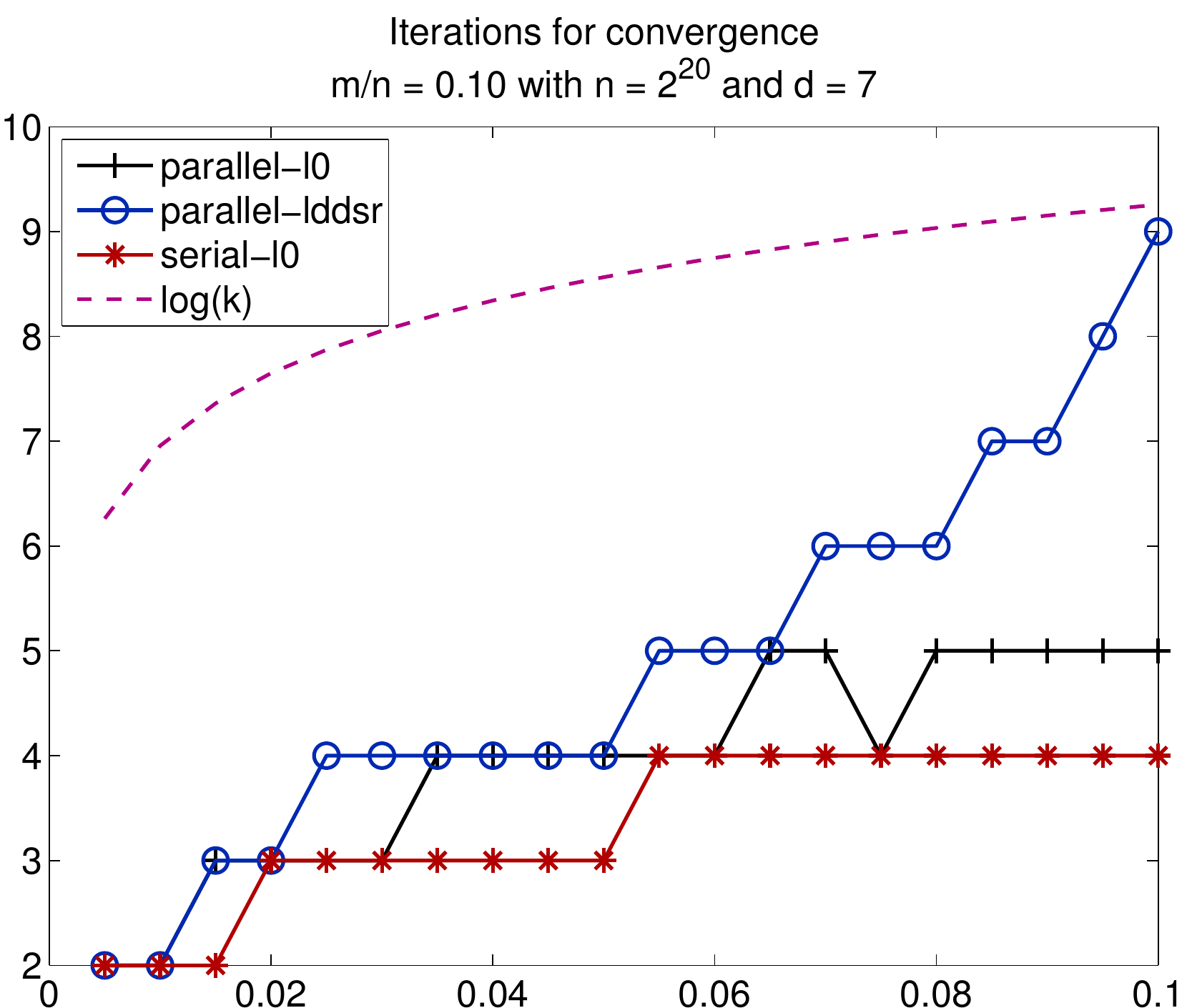}
	\caption[Iterations vs rho sublinear]{Number of iterations to convergence for Parallel-$\ell_0$, Serial-$\ell_0$, and parallel-LDDSR at $\delta = 0.1$ with $\mathbb{E}_{k, \varepsilon, 7}$ and $n = 2^{20}$. The number of iterations of Serial-$\ell_0$ has been normalised by $n$ to showcase its $\mathcal{O}(n\log k)$ guarantee in the number of iterations.}
	\label{fig:iterations-sublinear}
\end{figure}

\begin{figure}[!htbp]
	\centering
	\includegraphics[scale=0.39]{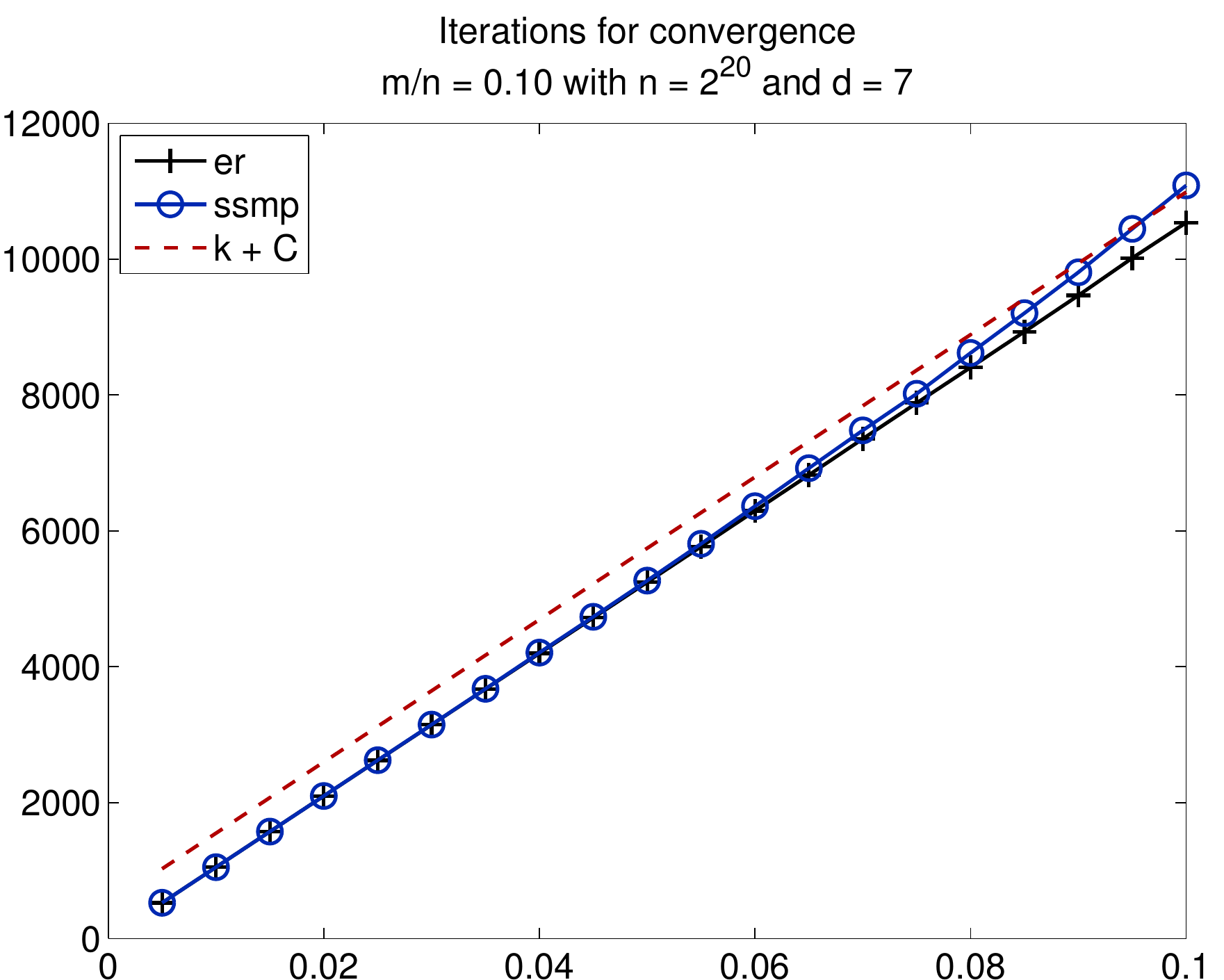}
	\caption[Iterations vs rho linear]{Number of iterations to convergence for ER and SSMP at $\delta = 0.1$ with $\mathbb{E}_{k, \varepsilon, 7}$ and $n = 2^{20}$.}
	\label{fig:iterations-linear}
\end{figure}

The theoretical guarantees of Serial-$\ell_0$ and Parallel-$\ell_0$
state that convergence can be achieved in $\mathcal{O}(nd \log k)$ operations.
The number of operations per iteration can be verified simply by
counting operations in the algorithm, which is $\mathcal{O}(d)$ for
Serial-$\ell_0$ and $\mathcal{O}(nd)$ for Parallel-$\ell_0$ and
recording the number of iterations.
Figure \ref{fig:iterations-sublinear} shows that the number of iterations to convergence for Serial-$\ell_0$, Parallel-$\ell_0$, and parallel-LDDSR. 
The tests were performed by fixing $n = 2^{20}$, $\delta = 0.1$, and $d = 7$, and considering signals with sparsity ranging from $\rho = 0.05$ to $\rho = 0.1$.
It can be seen in Figure \ref{fig:iterations-sublinear} that the number of iterations to convergence is bounded by the curve $f(k) = \log k$, thus verifying our claims.
We also make clear that by Definition \ref{def:iteration},
Serial-$\ell_0$ is shown to converge in $\mathcal{O}(n\log k)$
iterations, but for the sake of this experiment, we normalise the
final number of iterations for Serial-$\ell_0$ by a factor of $n$.  Note the lower number
of iteration by Serial-$\ell_0$ due to its serial implementation with
residual updates revealing more entries that satisfy the reduction of
the residual by $\alpha$.
Now, to give a point of comparison, we also compute the number of iterations for ER and SSMP, which take $\mathcal{O}(k)$ iterations to converge.
The results are shown in Figure \ref{fig:iterations-linear}, where the same parameters as in Figure \ref{fig:iterations-sublinear} have been used.
In particular, we can see that for a problem with $k/m = 0.8$, Parallel-$\ell_0$ takes 5 iterations, while ER and SSMP take about 8000 iterations to solve the same problem.

\subsection{Increasing phase transition as $\delta \rightarrow 0$ and $n \rightarrow \infty$}

\begin{figure}[!htbp]
\centering
  \includegraphics[scale=0.39]{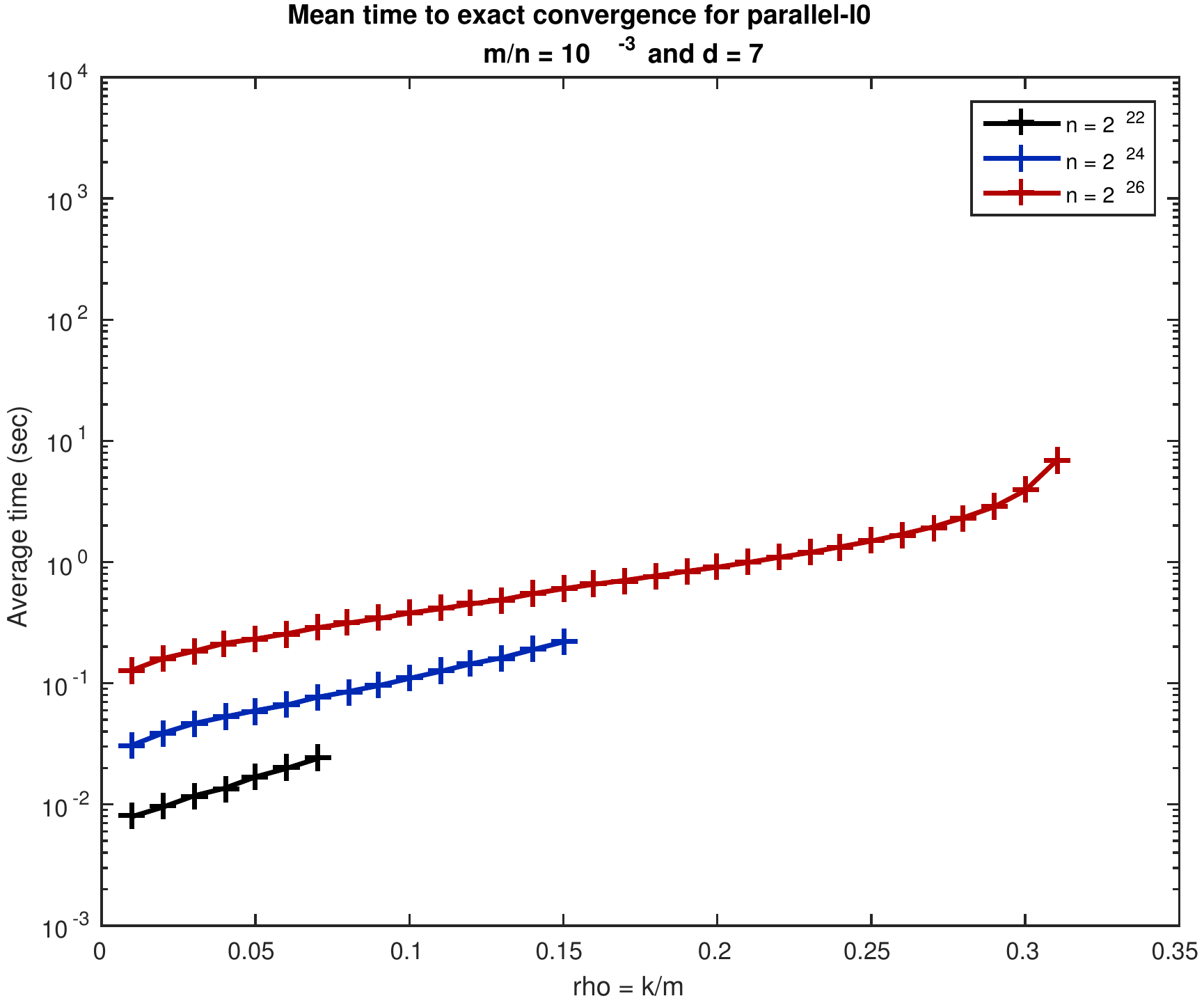}
   \caption{Average recovery time (sec) for Parallel-$\ell_0$, with dependence on $\rho$ for $\delta = 0.001$ and $\mathbb{E}_{k, \varepsilon, 7}$ with $n \in \{2^{22}, 2^{24}, 2^{26}\}$.}
    \label{fig:timing_low_delta}
\end{figure}

It is shown in Figure \ref{fig:phase-transitions} that Serial-$\ell_0$
and Parallel-$\ell_0$ have a very high phase transition of just over
0.3 even for very small values of $\delta$.
We hypothesise that this high phase transition persists for any fixed
$\delta\in (0,1)$ provided $n$
is sufficiently large.  
We provide numerical support of this claim in Figure \ref{fig:timing_low_delta}, where for fixed $\delta = 10^{-3}$ and $d = 7$, we have plotted the average time to convergence for Parallel-$\ell_0$ as $\rho$ increases.
The experiment was repeated for each $n \in \{2^{22}, 2^{24}, 2^{26}\}$, by initialising $\rho = 0.01$ and generating 30 problems at each $\rho$.
If at least 50\% of the problems converge we average out the time to convergence for successful cases, and perform the update $\rho \leftarrow \rho + 0.01$; otherwise, we stop.
Our results in Figure \ref{fig:timing_low_delta} show that for $\delta
= 10^{-3}$, the phase transition of the algorithm increases with $n$
to just over $0.3$.

Finally, in Table \ref{table:time_low_delta} we show the average
timing depicted in Figure \ref{fig:timing_low_delta} for $\rho = 0.05$
which shows the approximate increase in the average computation time
being proportional to $n$.
\begin{table}[!htbp]
\begin{center}
\begin{tabular}{|c|c|c|}
\hline
$n$& time $t_{n}$ & ratio $t_{4n}/t_{n}$\\
\hline
&&\\
$2^{22}$ & 0.0167 & 3.338\\
&&\\
$2^{24}$ & 0.0557 & 4.163\\
&&\\
$2^{26}$ & 0.2319 & - \\
&&\\
\hline
\end{tabular}
\end{center}
\caption{Average recovery time (sec) for Parallel-$\ell_0$ at $\rho = 0.05$ and $\delta = 10^{-3}$ for $n \in \{2^{22}, 2^{24}, 2^{26}\}$.}
\label{table:time_low_delta}
\end{table}

\subsection{Almost dissociated signals}

\begin{figure}[!htbp]
	\centering
	\includegraphics[scale=0.39]{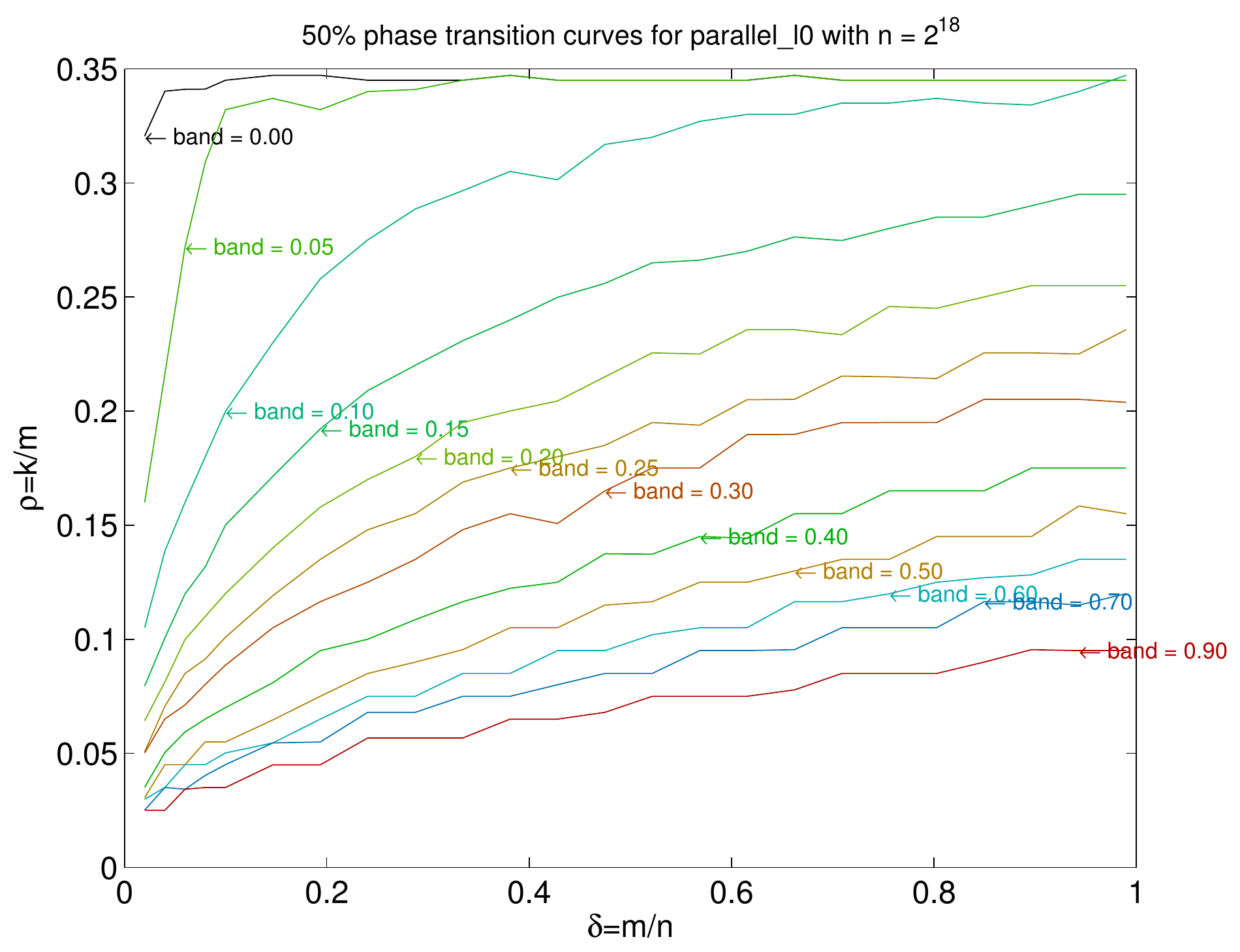}
	\caption[Parallel-$\ell_0$ several $d$'s]{50\% recovery probability logistic regression curves for Parallel-$\ell_0$ with $\mathbb{E}_{k, \varepsilon, 7}$ and $n = 2^{18}$, with signals having a fixed proportion, {\em band}, of identical nonzero elements in its support.}
	\label{fig:banded_signals}
\end{figure}

The analysis of Parallel-$\ell_0$ and Serial-$\ell_0$ relied on the
model of dissociated signals \eqref{eq:dissociated}.  
We explore the effect on recovery ability of Parallel-$\ell_0$ and
Serial-$\ell_0$ as the signal model is no longer dissociated, with a
fixed fraction of the values in $x$ being equal.
To do this, we consider signals $x \in \Sparse$ with nonzero values composed of two bands: one in which {\em all} entries are equal to a fixed value drawn at random from a standard normal distribution $\mathcal{N}(0, 1)$, and another one in which {\em each} entry is drawn independently of each other from $\mathcal{N}(0,1)$.
Our results are shown in Figure \ref{fig:banded_signals}, where we can
see that as the fraction of values which are equal increases (shown in the
figure by the parameter {\em band}), the phase transitions gracefully
decrease from the flat shape observed for perfectly dissociated
signals to an increasing log-shaped curve when $\mbox{\em band }=
0.9$.  Note that the overall phase transition decreases, with the
greatest decrease for $\delta\ll 1$.

\subsection{$d$ should be small, but not too small}

\begin{figure}[!htbp]
	\centering
	\includegraphics[scale=0.39]{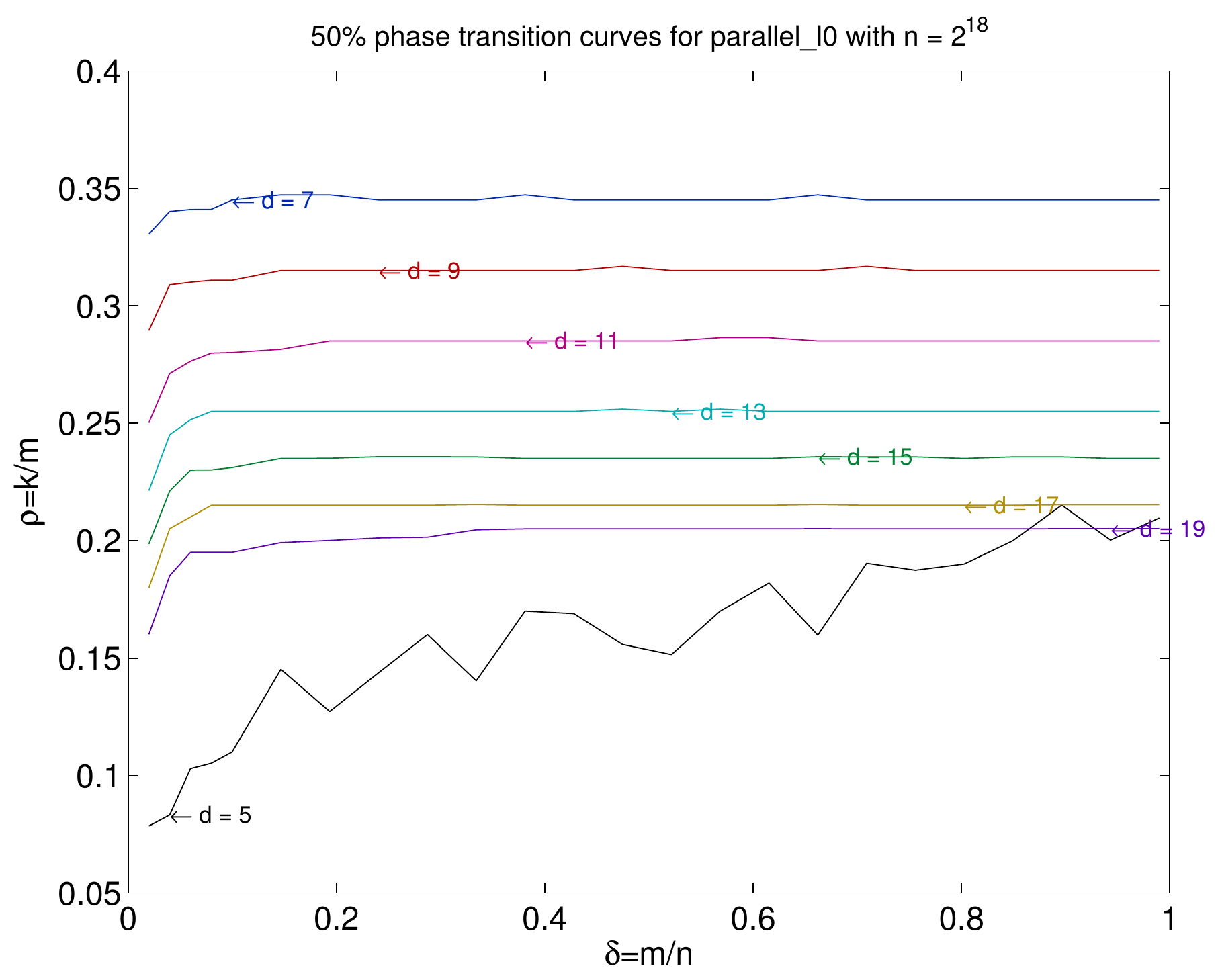}
	\caption[Parallel-$\ell_0$ several $d$'s]{50\% recovery probability logistic regression curves for Parallel-$\ell_0$ with $\mathbb{E}_{k, \varepsilon, d}$ and $n = 2^{18}$ for $d \in \{5, 7, 9, 11, 13, 15, 17, 19\}$.}
	\label{fig:parallel_l0_several_d}
\end{figure}

Selection of the number of nonzeros per column, $d$, has not been adressed.
In our numerical experiments we have consistently chosen $d = 7$ as
the left-degree of our expander.
Our choice of $d=7$ for our problem size's order of magnitude is
justified by Figure \ref{fig:parallel_l0_several_d}, where we have
computed the phase transitions for Parallel-$\ell_0$ for all odd
values of $d$ between 5 and 19.
For $d = 5$, the phase transition of the algorithm is very low, thus
signalling expanders of bad quality. For $d=7$ the phase transition
is substantially greater than when $d=5$, and gradually decreases for values of $d$
greater than seven.   Note that the expander condition implies
$(1-\varepsilon)dk<m$ which encourages small values of $d$ in order
that $m/k$ can be as large as possible.

\section{Conclusions and Future Work}
\label{sec:conclusions}

We have proposed two algorithms for combinatorial compressed sensing with provable convergence guarantees in $\mathcal{O}(dn \log k)$ operations and very high phase transitions when the signal $x$ is dissociated. In particular, Parallel-$\ell_0$ is observed to be empirically the fastest algorithm in compressed sensing when the signal is dissociated. We have used the dissociated signal model in the convergence proofs, but that in practice one can relax this assumption and still get reasonably high phase transitions.  

As future work it remains to address the case of noisy observations, and to extend the scope of the algorithms to more general signal models.  The proofs presented in this paper should extend trivially to noise which is bounded to be less than half the minimal distance between obtainable values $\sum_{i\in T} x_i$ by introducing an equivalence class.  A variant which is robust to Gaussian noise is scope for future work.

\appendix
\label{app:expander-proofs}

For completeness, we give a proof of Lemma \ref{th:equivalence-lemma}

\begin{IEEEproof}
For any unbalanced, left $d$-regular, bipartite graph it holds that:
\begin{equation}
\label{eq:partition}
|\mathcal{N}_1(S)| + |\mathcal{N}_{>1}(S)| = |\mathcal{N}(S)|,
\end{equation}
\begin{equation}
\label{eq:double-counting}
|\mathcal{N}_1(S)| + 2|\mathcal{N}_{>1}(S)| \leq d|S|.
\end{equation}
Where (\ref{eq:partition}) follows from the definition of $\mathcal{N}_{>1}(S)$, and (\ref{eq:double-counting}) by double-counting the edges emanating from $S$ to $\mathcal{N}(S)$. Now, to prove that (\ref{eq:unique-neighbors}) is necessary, assume that $G$ is a $(k, \varepsilon, d)$-expander graph. Then, for $S \in [n]^{(\leq k)}$ we have that
\begin{equation}
\label{eq:expander-1}
|\mathcal{N}(S)| > (1 - \varepsilon)d|S|.
\end{equation}
Combining (\ref{eq:partition}), (\ref{eq:double-counting}) and (\ref{eq:expander-1}) we get the chain of inequalities
\begin{equation}
d|S| - |\mathcal{N}_{>1}(S)| \geq |\mathcal{N}(S)| > (1 - \varepsilon)d|S|,
\end{equation}
which yield
\begin{equation}
\label{eq:intermediate-1}
|\mathcal{N}_{>1}(S)| < \varepsilon d|S|.
\end{equation}
Plugging (\ref{eq:intermediate-1}) into (\ref{eq:partition}) and using (\ref{eq:expander-1}) we obtain
\begin{equation}
\label{eq:expander-2}
|\mathcal{N}_1(S)| > (1 - 2 \varepsilon)d|S|.
\end{equation}
To prove the sufficiency of (\ref{eq:unique-neighbors}) for graph expansion, we couple it with (\ref{eq:double-counting}) into the system
\begin{equation}
(1 - 2\varepsilon)d|S| < |\mathcal{N}_1(S)| \leq d|S| - 2|\mathcal{N}_{>1}(S)|,
\end{equation}
and use the left and right hand sides recover (\ref{eq:intermediate-1}). Now, using (\ref{eq:unique-neighbors}) and (\ref{eq:partition}) we obtain
\begin{equation}
\label{eq:intermediate-2}
|\mathcal{N}(S)| - |\mathcal{N}_{>1}(S)| > (1 - 2\varepsilon)d|S|.
\end{equation}
And using (\ref{eq:intermediate-1}) in (\ref{eq:intermediate-2}) allows us to recover (\ref{eq:expander-1}), implying that $G$ is a $(k, \varepsilon, d)$-expander graph.
\end{IEEEproof}

%


\bibliographystyle{IEEEtran} 
\bibliography{expander_parl0}

\end{document}